\documentclass[final,5p,12pt,fleqn]{elsarticle}
\makeatletter
\def\ps@pprintTitle{%
	\let\@oddhead\@empty
	\let\@evenhead\@empty
	\def\@oddfoot{}%
	\let\@evenfoot\@oddfoot
}
\makeatother
\usepackage{amssymb}
\usepackage{amsmath}
\usepackage{placeins} 	%MS package needed for the float barrier
\usepackage{url} 		%MS package needed to break urls in the references

\usepackage{tabularx}   %MS package needed for tables
\usepackage{booktabs}   %MS package for toprule midrule etc.
\usepackage{todonotes}  %MS package for todo notes
\usepackage{siunitx}    %MS package for units 
\usepackage{tablefootnote} %MS package for table footnotes
\usepackage{makecell}     % MS package for merged headers/multicells in tables
\usepackage{eurosym}
\DeclareSIUnit{\sieuro}{\mbox{\euro}}
\DeclareSIUnit{\ct}{ct}
\DeclareSIUnit\year{yr}

\usepackage{caption} % MS for separate captions within mini pages
\usepackage[export]{adjustbox} % MS needed for figure horizontally next to table magic..
\usepackage{bm} % MS required to write bold in math mode 

\usepackage{tikz} %MS
\usetikzlibrary{arrows.meta, positioning, shapes.geometric, fit, calc} %MS

\journal{Elsevier}

\begin{document}
\emergencystretch 3em %MS used to avoid overfull lines (lines that are going beyond the columnwidth)

\begin{frontmatter}


\title{Optimal Heat Storage Sizing for District Heating Networks to Maximize Electricity Revenue from Combined Heat and Power Units}

\author[1,3]{Martin Sollich\fnref{label1}}
\author[1,3]{Maarten Blommaert}% \fnref{label1,label3}}

\affiliation[1]{organization={Department of Mechanical Engineering, KU Leuven},addressline={ Celestijnenlaan 300 box 2421},city={Leuven},postcode={3001},country={Belgium}}

\affiliation[3]{organization={EnergyVille},addressline={Thor Park, Poort Genk 8310},city={Genk},postcode={3600},country={Belgium}}

\fntext[label1]{Corresponding author. Email address: martin.sollich@kuleuven.be}

%% Abstract
\begin{abstract}
Integrating heat storages in district heating networks (DHNs) supports managing dynamic characteristics such as heat-demand variations and changing energy prices, and the intermittency of renewable sources. A common application are DHNs with combined heat and power (CHP) units, where a storage allows shifting heat extraction to periods with favorable electricity prices. Although the benefits of heat storage in DHNs are well established, determining the optimal storage size remains challenging due to the diversity of DHNs in terms of heat sources, production and consumption patterns, network heat losses, and fuel costs. This paper presents a scalable, automated methodology for optimizing short-term heat storage in DHNs using mathematical optimization. The nonlinear, physics-based approach models the DHN, heat producers, and storages simultaneously to minimize total economic cost through optimal storage sizing, explicitly considering time-varying heat-production costs, heat demands, and heat losses. The methodology is demonstrated on a 3rd-generation CHP-DHN in Belgium with 30 consumers, exceeding the scale of previous physics-based storage-sizing studies. For this system, optimal storage integration reduces the 20-year cost by $730\si{\,\kilo\sieuro}$ (16.5\%), from $4.41\si{\,\mega\sieuro}$ to $3.68\si{\,\mega\sieuro}$. The reduction results from a $1.07\si{\,\mega\sieuro}$ decrease in heat-production cost achieved by shifting heat extraction to periods of low electricity prices, while the storage investment amounts to $323\si{\,\kilo\sieuro}$. A comparison to a commonly used simplified storage sizing method, which fails to identify the optimal storage size, demonstrates the advantage of the proposed holistic, physics-based optimization approach in ensuring feasible designs, accurate cost assessments, and optimal storage sizing.
\end{abstract}

\begin{graphicalabstract}
%twocolumnfigure
\includegraphics[width=2\columnwidth]{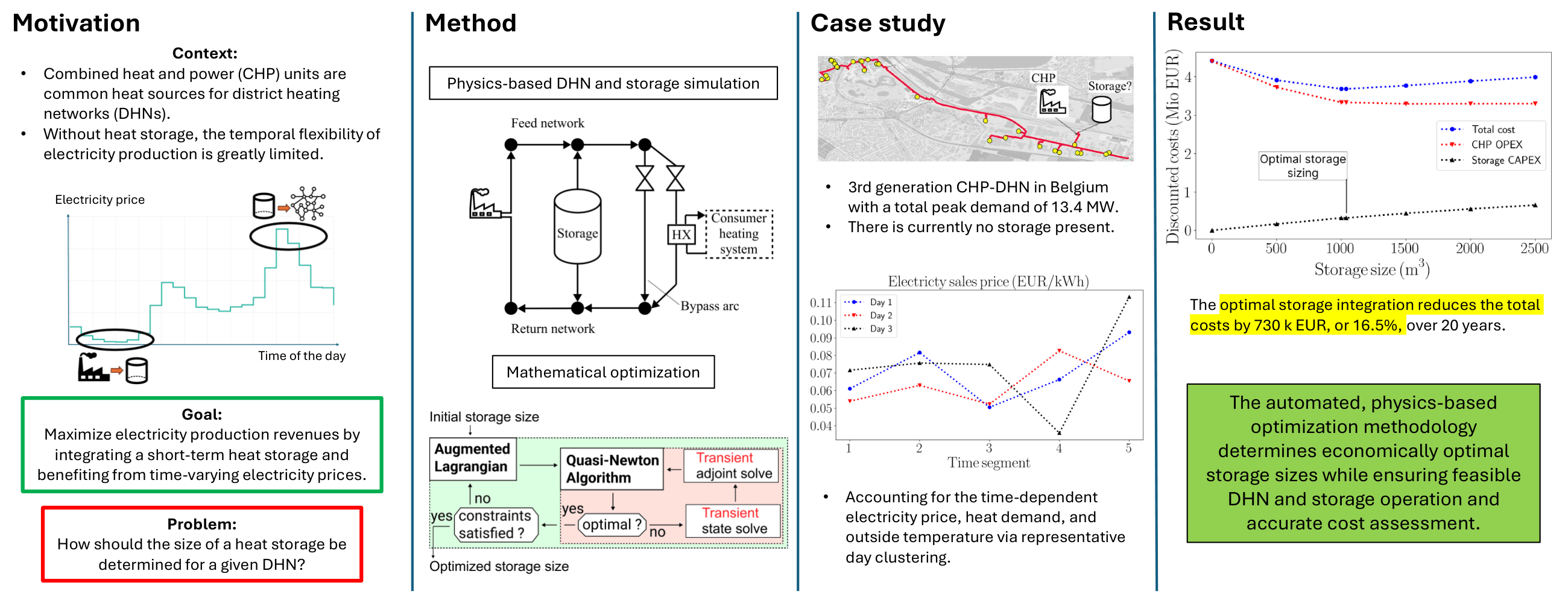}
\end{graphicalabstract}

\begin{highlights}
\item Physics-based heat storage optimization for district heating networks is proposed.
\item Nonlinear mathematical optimization automates heat storage sizing.
\item Storage reduces total costs by shifting heat extraction to moments of low prices.
\item Method is demonstrated on a real combined heat and power district heating network.
\item The optimal storage integration reduces the 20-year cost by $730\si{\,\kilo\sieuro}$ or 16.5\%.
\end{highlights}

%% Keywords

\end{frontmatter}

\FloatBarrier


\section{Introduction}
\label{sec:introduction}

%add a table only with the abbreviations (not with the variables, subscripts etc.)
\begin{table}[t]
	\centering
	\caption{Abbreviations used in this manuscript.}
	\label{tab:abbreviations}
	\begin{tabularx}{\columnwidth}{>{\hsize=.34\hsize}X>{\hsize=.66\hsize}X}
		\toprule
		\textbf{Abbreviation} & \textbf{Description} \\
		\midrule
		CAPEX & Capital expenditures \\
		CHP & Combined heat and power \\
		CHP-DHN & CHP district heating network\\
		DH & District heating \\
		DHN(s) & District heating network(s) \\
		GIS & Geographic information system \\
		HX & Heat exchanger \\
		HS & Heat storage \\
		LCOH & Levelized cost of heating \\
		LP & Linear program \\
		MILP & Mixed-integer linear program \\
		NLP & Nonlinear program \\
		OPEX & Operating expenditures \\
		\bottomrule
	\end{tabularx}
\end{table}

\newcommand{\COTWO}{\mathrm{CO_2}}

The European Commission identifies district heating (DH) as a cost-effective option to reduce fossil fuel dependence in the heating sector \cite{REPowerEU}. In district heating networks (DHNs), heat storages enable improved handling of dynamic system characteristics, including variations in heat and electricity demand, energy prices, renewable intermittency, and system disturbances \cite{GUELPA2019}.

A common application of heat storage in DHNs is in combination with combined heat and power (CHP) units, which are widely used as heat sources. Storage integration increases operational flexibility, allowing higher electricity sales revenues and reducing required CHP and peak boiler capacity through peak shaving \cite{GUELPA2019}.

Despite these benefits, determining the optimal storage size remains challenging due to differences between DHNs in heat sources, seasonal production and consumption patterns, and fuel costs \cite{SIFNAIOS2025}. Consequently, optimization-based methods have been proposed to identify cost-optimal storage sizes \cite{Wang2015,Morvaj2016}.
 
\subsection{Previous work on optimal storage sizing for DHNs}
\label{subsec:previous work on storage sizing}

Optimization approaches for heat storage sizing in DHNs can be grouped into two main categories. The first linearizes the nonlinear heat transport problem and cost functions (simplified modeling approaches), while the second accounts for non-linearities and explicitly represents the underlying physics (physics-based nonlinear modeling approaches).

The following literature review focuses on optimal sizing of short-term heat storage, which is the main scope of this work. Some studies on seasonal storage are also included, as research on short-term storage sizing remains limited and several works address both short-term and seasonal storage.

\subsubsection{Simplified modeling approaches}
\label{subsec:previous work storage sizing - simplified} 

\citet{Wang2015} solved a linear program (LP) to optimize CHP and solar thermal operation with heat storage, considering two fixed storage sizes and both short-term and seasonal storage, while aggregating the DHN into a single demand point. \citet{Morvaj2016} proposed a multi-objective mixed-integer linear program (MILP) including economic and environmental criteria for a small network of 11 houses, optimizing production technologies, storage capacity, network layout, and operation using representative design days. \citet{VANDERHEIJDE2019} introduced a MILP-based heuristic method to select the storage size from a predefined discrete set, focusing on seasonal storage in a small-scale network with 3 demand points. \citet{JEBAMALAI2020} developed a linear, heuristic approach to study the effect of storage size and location in a DHN with 7 substations, with storage sizing performed parametrically. \citet{Resimont2021} presented a multi-period MILP for DHN design that jointly optimized storage sizing and location, production units, and network layout, demonstrated on a small network with 16 edges.
 
\citet{PANS2024} used a simulation-based Excel model to compare predefined scenarios with and without storage, analyzing decentralized short-term storage and seasonal storage impacts on costs and $\COTWO$ emissions, while the seasonal storage size was varied parametrically. More recently, \citet{CERUTI2025} proposed a spatially resolved MILP, optimizing network layout, investments, and DHN operation, including short-term storage, and applied it to a DHN with 122 consumers.

These linear optimization approaches neglect the inherently nonlinear physics of heat and mass transfer in DHNs and storage systems, as well as nonlinear cost functions such as storage investment costs. As a result, they cannot guarantee feasibility or provide accurate cost assessments for DHN and heat storage design and operation.

\subsubsection{Physics-based nonlinear modeling approaches}
\label{subsec:previous work storage sizing - nonlinear}

Among nonlinear approaches, several studies apply simulation-based parameter optimization. \citet{SALOUX2021} conducted a study for a solar-based DHN aggregated into a single demand point, optimizing solar and storage operation for energy savings while investigating short-term and seasonal storage capacities. \citet{BELLOS2022} analyzed decentralized short-term storage in a small DHN with 9 buildings, comparing central and decentralized configurations for domestic hot water preparation over a single-day horizon. \citet{SIFNAIOS2025} proposed a simulation-based economic parameter optimization including short-term and seasonal storage, optimizing the capacities of storage, heat pumps, and electric boilers. They used one- to three-week control horizons combined with full-year simulations and aggregated the DHN into a single demand point.

In addition, automated optimization approaches based on heuristic or mathematical optimization have been proposed. \citet{TALEBI2019} applied a simulation-based heuristic method to size producers, short-term and seasonal storage for an aggregated (i.e. single demand point) DHN. \citet{ABOKERSH2020} developed a nonlinear heuristic optimization for a small DHN with 10 houses, optimizing seasonal and domestic hot water storage volumes. \citet{MAXIMOV2021} used a genetic algorithm to optimize solar thermal and storage capacities with cost and emission objectives, considering no more than 6 design variables. \citet{Delubac2023} proposed a heuristic framework to optimize solar thermal systems, short-term storage, and peak boilers for an aggregated DHN using representative days. \citet{Velasco2023} compared predefined producer configurations with and without storage in terms of costs and $\COTWO$ emissions, while \citet{HASSAN2024} presented a nonlinear mathematical optimization that optimized producer and short-term storage sizes and locations using representative days, with the DHN being aggregated into a single demand point in both studies.

No existing study demonstrates an automated, physics-based storage sizing approach that simultaneously models nonlinear physics of heat transfer in the DHN and heat storage, while applying scalable mathematical optimization. Existing nonlinear methods either aggregate the DHN into a single demand point or are limited to small networks with up to 10 consumers. To find the optimal storage size, the feasibility of satisfying heat demand and meeting network pressure limits must be ensured by modeling nonlinear heat and pressure losses throughout the DHN. This modeling must account for existing pipe infrastructure and the nonlinear physics of the storage. Moreover, the storage size and operation must be determined automatically and in a scalable way through mathematical optimization, taking into account nonlinear cost functions.

\subsection{Goal and scope of this paper}
\label{subsec:goal and scope}

This paper presents a scalable, automated optimization methodology for cost-optimal integration of short-term heat storage in DHNs. The nonlinear, physics-based approach simultaneously models the DHN, heat producers, and storages in detail to minimize operator costs by optimally sizing storage, while accounting for time-varying heat production costs and heat demand. The methodology is demonstrated on a real 3rd-generation CHP-DHN in Belgium, where a central short-term storage is integrated to exploit time-varying electricity prices.

Existing storage sizing approaches either simplify the storage and network physics entirely or model the storage physics accurately in a nonlinear fashion while aggregating the entire network into a single demand point. This can lead to sub-optimal storage sizing when network heat losses are not accurately accounted for, or may require time-consuming iterations between isolated simulation and optimization steps. Moreover, the absence of spatially-resolved physical network information inherently limits applicability, precluding storage integration studies in which potential storage locations throughout the network are to be assessed. 

The novelty of this work lies in the automated, physics-based optimization of short-term heat storage sizing for DHNs. In contrast to existing approaches, the full spatial detail of the network physics -- momentum and heat losses -- is modeled jointly with producers and storage in a nonlinear formulation, ensuring feasibility, accurate cost evaluation, and optimal storage sizing while remaining scalable to large networks. Scalability is achieved through three key modeling choices: transforming the intrinsically discrete charging and discharging decision into a continuous problem formulation; modeling transient behavior only where necessary -- for the storage -- while retaining a quasi-steady approach for the DHN; and computing gradients efficiently via the adjoint method. Furthermore, this work is the first to apply an automated, physics-based storage optimization method to a CHP-DHN to exploit time-varying electricity prices. It builds on the economic DHN topology optimization framework first introduced by \citet{Blommaert}, which has since been extended to enable, among others, optimal producer retrofits for reducing network temperatures and $\COTWO$ emissions \cite{SOLLICH2025Producer}. The current paper aims to extend the approach to include heat storage by developing a storage model, the required transient simulation, and the required transient adjoint-gradient calculation.

\newcommand{\ve}[1]{\bm{#1}} 	% vector (bold font)
\newcommand{\tp}[1]{#1^{\intercal}} 			% transpose

\newcommand{\gEdge}[3]{#1_{#2#3}}   % generic edge
\newcommand{\gNode}[2]{#1_{#2}}     % generic node

% Directed graph definition
\newcommand{\dirGraph}{G}
\newcommand{\setNodes}{N}
\newcommand{\setEdges}{E}

% Graph components in district heating networks
\newcommand{\pro}{\mathrm{pr}}
\newcommand{\con}{\mathrm{con}}
\newcommand{\stor}{\mathrm{stor}}
\newcommand{\hx}{\mathrm{hx}}
\newcommand{\hp}{\mathrm{hp}}
\newcommand{\byp}{\mathrm{bp}}
\newcommand{\jun}{\mathrm{jun}}
\newcommand{\pipe}{\mathrm{pipe}}
\newcommand{\operation}{\mathrm{op}}
\newcommand{\storHot}{\mathrm{h}}
\newcommand{\storCold}{\mathrm{c}}
\newcommand{\totalMassIDX}{\mathrm{tot}}

%% subsets
% nodes
\newcommand{\Npro}{\setNodes_\pro}
\newcommand{\Ncon}{\setNodes_\con}
\newcommand{\NconF}{\setNodes_{\con,\mathrm{f}}}
\newcommand{\NconR}{\setNodes_{\con,\mathrm{r}}}
\newcommand{\Njun}{\setNodes_\jun}
\newcommand{\NproF}{\setNodes_{\pro,\mathrm{f}}}
\newcommand{\NproR}{\setNodes_{\pro,\mathrm{r}}}
\newcommand{\Nstor}{\setNodes_\stor}
\newcommand{\NstorF}{\setNodes_{\stor,\mathrm{f}}}
\newcommand{\NstorR}{\setNodes_{\stor,\mathrm{r}}}

% edges
\newcommand{\EF}{\setEdges_{\mathrm{f}}}
\newcommand{\Epro}{\setEdges_\pro}
\newcommand{\EproGB}{\setEdges_{\pro,\GBIndex}}
\newcommand{\EproHP}{\setEdges_{\pro,\hp}}
\newcommand{\EproEB}{\setEdges_{\pro,\EBIndex}}
\newcommand{\EproST}{\setEdges_{\pro,\STIndex}}
\newcommand{\CHPIndex}{\mathrm{chp}}
\newcommand{\EproCHP}{\setEdges_{\pro,\CHPIndex}}

\newcommand{\Econ}{\setEdges_\con}
\newcommand{\Econhx}{\setEdges_{\con,\hx}}
\newcommand{\Econhp}{\setEdges_{\con,\hp}}
\newcommand{\Econbyp}{\setEdges_{\con,\byp}}
\newcommand{\Epipe}{\setEdges_\pipe}
\newcommand{\EpipeF}{\setEdges_{\pipe,\mathrm{f}}}
\newcommand{\EpipeR}{\setEdges_{\pipe,\mathrm{r}}}
\newcommand{\Eop}{\setEdges_\operation}

\newcommand{\Estor}{\setEdges_\stor}
\newcommand{\EstorHot}{\setEdges_{\stor,\storHot}}
\newcommand{\EstorCold}{\setEdges_{\stor,\storCold}}

% i,j definition
\newcommand{\gi}{i}
\newcommand{\gj}{j}

\newcommand{\giNode}[1]{\gNode{#1}{\gi}}
\newcommand{\gjNode}[1]{\gNode{#1}{\gj}}
\newcommand{\gijEdge}[1]{\gEdge{#1}{\gi}{\gj}}

\newcommand{\gjNodetime}[1]{\gNode{#1}{\gj,\timeVar}}
\newcommand{\giNodetime}[1]{\gNode{#1}{\gi,\timeVar}}
\newcommand{\gijEdgetime}[1]{\gEdge{#1}{\gi}{\gj},_\timeVar}

%% NOTATION %%
%% Optimization problem
\newcommand{\cost}{J}
\newcommand{\equalCon}{\ve{c}}
\newcommand{\inEqualConSkalar}{\ve{h}}
\newcommand{\inEqualConSkal}{h}
\newcommand{\EqualConSkal}{g}
\newcommand{\inEqualCon}{\inEqualConSkalar}
%%% Variables
%% design variables
\newcommand{\designVarskal}{\varphi}
\newcommand{\designVar}{\ve{\designVarskal}}
\newcommand{\designVarTime}{\ve{\varphi}_{\timeVar}}
\newcommand{\designVarDT}{\ve{\varphi}_{\dayTimeVar}}
\newcommand{\designVarTimeInvariantSkalar}{\phi}
\newcommand{\designVarTimeInvariant}{\ve{\designVarTimeInvariantSkalar}}
\newcommand{\designVarTimeDependentSkalar}{\bar{\phi}}
\newcommand{\designVarTimeDependent}{\ve{\designVarTimeDependentSkalar}}
%% state variables
\newcommand{\stateVarSkal}{x}
\newcommand{\stateVar}{\ve{\stateVarSkal}}
\newcommand{\stateVarTime}{\ve{\stateVarSkal}_{\timeVar}}
\newcommand{\stateVarDT}{\ve{\stateVarSkal}_{\dayTimeVar}}
\newcommand{\topVarSkalar}{d}
\newcommand{\topVar}{\ve{\topVarSkalar}}
% physical state variables
\newcommand{\stateVarTemp}{\theta}
\newcommand{\temperature}{T}
\newcommand{\stateVarPress}{p}
\newcommand{\stateVarFlow}{q}
\newcommand{\stateVarHeat}{\dot{Q}}
\newcommand{\stateVarStorMass}{m}
\newcommand{\defFlow}{\ve{\stateVarFlow}}
\newcommand{\defPressure}{\ve{\stateVarPress}}
\newcommand{\defTemp}{\ve{\stateVarTemp}}
\newcommand{\defStorMass}{\ve{\stateVarStorMass}}

\newcommand{\defStateVar}{\ve{\stateVar} = \tp{\left[\ve{\stateVarStorMass},\ve{\stateVarPress},\ve{\stateVarFlow},\ve{\stateVarTemp}\right]}}
%design variables
\newcommand{\capVar}{\phi} 	% design variable producer inflow
\newcommand{\prodInput}{\delta} 	% design variable producer inflow
\newcommand{\prodTemp}{\tau}     % producer temp design var
\newcommand{\HXValve}{\alpha}    % HX valve design var
\newcommand{\bypValve}{\beta}    % bypass valve design var
\newcommand{\HPValve}{\gamma}    % HP valve design var
\newcommand{\subSelection}{\psi} % substation selection
\newcommand{\subSelectionPenalized}{\bar{\psi}} % substation selection
\newcommand{\storSize}{\zeta} % storage size
\newcommand{\storFlow}{\iota} % storage charging flow

% storage parameters
\newcommand{\storMaxVolume}{V_{\stor,\mathrm{max}}} % storage charging flow
\newcommand{\Ustor}{U_{\stor}}
\newcommand{\storSurface}{A}

% time periods, days, segments etc.
\newcommand{\timeVar}{t} % this is the time variable, used when writing time-dependent differential equations with (t), and NOT the time segment index
\newcommand{\timeVarPlusOne}{t+1}
\newcommand{\deltaTime}{\Delta\timeVar}
\newcommand{\initINDEX}{\mathrm{init}}
\newcommand{\timeVarPeak}{t_{\mathrm{peak}}}
\newcommand{\setPeriods}{\Upsilon} % used for temporal clustering methods other than representative day clustering where the periods are fully independent of each other and simply defined via one set
\newcommand{\peakPeriod}{t_{\textrm{peak}}}

\newcommand{\dayVar}{d} % indicating the representative day in case of representative day clustering
\newcommand{\timeSegVar}{s} % indicating the time segment in case of representative day clustering
\newcommand{\dayTimeVar}{\dayVar,\timeSegVar} % for easier access

% set for the time segments and day in case of representative day clustering
\newcommand{\timeVarSetd}{\mathcal{S}_{d}} % set of time segments dependent on the day d 
\newcommand{\dayVarSet}{\mathcal{D}} 

\newcommand{\n}{n}
\newcommand{\npipes}{\n_{\pipe}}
\newcommand{\nperiods}{\n_{\textrm{period}}} %number of periods
\newcommand{\ndays}{\n_{\textrm{day}}}    %number of days (days are periods are the same)
\newcommand{\nts}{\n_{\textrm{ts}}} % number of time segments per day 
\newcommand{\ntsd}{\n_{\textrm{ts},d}} % number of time segments per day dependent on the day d
\newcommand{\timeSlices}{\nperiods}

\newcommand{\defDesignVar}{\designVar=\tp{\left[\ve{\HXValve},\ve{\bypValve},\ve{\prodInput},\ve{\storSize},\ve{\storFlow}\right]}}

% Constraints
\newcommand{\defModelConstraints}{\equalCon_\timeVar\left( \designVarTime,\stateVar_\timeVar\right)}
\newcommand{\defModelConstraintsDT}{\equalCon_{\dayTimeVar}\left( \designVarDT,\stateVar_{\dayTimeVar}\right)}
\newcommand{\designVarUpper}{\designVar_{\mathrm{up}}}
\newcommand{\designVarLower}{\designVar_{\mathrm{low}}}
\newcommand{\opVarBoxConstraints}[1]{\designVarLower\leq #1 \leq \designVarUpper}

\newcommand{\defnPeriods}{\nperiods \in \mathbb{N}_1}

\newcommand{\nElements}{\n_{\textrm{spat}}}
\newcommand{\techCon}{\ve{\inEqualConSkal}_{\mathrm{tech},\timeVar}}
\newcommand{\techConDT}{\ve{\inEqualConSkal}_{\mathrm{tech},\dayTimeVar}}
\newcommand{\techConEquality}{\ve{\EqualConSkal}_{\mathrm{tech},\timeVar}}
\newcommand{\defStateConstraints}{\techCon(\designVarTime,\stateVarTime)}
\newcommand{\defStateConstraintsDT}{\techConDT(\designVarDT,\stateVarDT)}
\newcommand{\defStateConstraintsEquality}{\techConEquality(\designVarTime,\stateVarTime)}

\newcommand{\nDiscretePipes}{\n_{\mathrm{D}}}
\newcommand{\setDefDiscreteDiameters}{\{\diameterDiscrete_{0},\dots,\diameterDiscrete_{\nDiscretePipes}\}}
\newcommand{\defnpipes}{\npipes = \card{\Epipe}}

%producer names
\newcommand{\GBName}{\mathrm{gas \ boiler}}
\newcommand{\HPName}{\mathrm{heat \ pump}}
\newcommand{\STName}{\mathrm{solar \ thermal}}
\newcommand{\EBName}{\mathrm{electric \ boiler}}
%italic for titles
\newcommand{\GBNameITALIC}{gas \ boiler}
\newcommand{\HPNameITALIC}{heat \ pump}
\newcommand{\STNameITALIC}{solar \ thermal}
\newcommand{\EBNameITALIC}{electric \ boiler}
%producer indices
\newcommand{\GBIndex}{\mathrm{GB}}
\newcommand{\HPIndex}{\mathrm{HP}}
\newcommand{\STIndex}{\mathrm{ST}}
\newcommand{\EBIndex}{\mathrm{EB}}

%% parameters
\newcommand{\deltaTHPsub}{\Delta T_{\mathrm{evap},\con,\hp}} % temperature drop along the DHN arc supplying the HP sub
\newcommand{\deltaTHPsubValue}{5}
\newcommand{\mindTsubHP}{\Delta T_{{\mathrm{min\,lift}},\con,\hp}} % minimal temperature lift between HP sub source and sink
\newcommand{\Tsource}{T_{\mathrm{source}}}
\newcommand{\Tsink}{T_{\mathrm{sink}}}
\newcommand{\mindTsubHPValue}{5}
\newcommand{\KOPEX}{K}
\newcommand{\KOPEXvalue}{8208} %h/year
%consumer side
\newcommand{\Tsechot}{\stateVarTemp_{2,\text{h}}}								
\newcommand{\Tseccold}{\stateVarTemp_{2,\text{c}}}
\newcommand{\TsechotIJ}{\stateVarTemp_{2,\text{h},ij}}								
\newcommand{\TseccoldIJ}{\stateVarTemp_{2,\text{c},ij}}
\newcommand{\TsechotdesignAbsolute}{T_{2,\text{h}}}	% secondary temperature hot
\newcommand{\TseccolddesignAbsolute}{T_{2,\text{c}}}	% secondary temperature cold

%consumer side

%network side
\newcommand{\Tphotdesign}{\dTinf_{1\text{,h}}}	% primary temperature hot
\newcommand{\Tpcolddesign}{\dTinf_{1\text{,c}}}	% primary temperature cold
\newcommand{\TphotdesignNom}{\dTinf_{1\text{,h,nom}}}	% primary temperature hot
\newcommand{\TpcolddesignNom}{\dTinf_{1\text{,c,nom}}}	% primary temperature cold
\newcommand{\TphotdesignAbsolute}{T_{1\text{,h}}}	% primary temperature hot
\newcommand{\TpcolddesignAbsolute}{T_{1\text{,c}}}	% primary temperature cold 
\newcommand{\TphotdesignAbsoluteNom}{T_{1\text{,h,nom}}}	% primary temperature hot
\newcommand{\TpcolddesignAbsoluteNom}{T_{1\text{,c,nom}}}	% primary temperature cold 

%% producer 
% temperature
\newcommand{\deltaTnetwork}{\Delta T_{\mathrm{primary}}}
% flow rate
\newcommand{\qnetwork}{\mathrm{Network flowrate}}

\newcommand{\Cmin}{C_{\mathrm{min}}}										% minimum heat capacity
\newcommand{\Cmax}{C_{\mathrm{max}}}										% maximum heat capacity
\newcommand{\Cstar}{C^*}										% ratio of heat capacities
\newcommand{\NTU}{NTU}											% number of transfer units
\newcommand{\U}{U}												% overall heat transfer coefficient
\newcommand{\epsilonNTU}{\epsilon}	
\newcommand{\Qdemand}{\dot{Q}_{\mathrm{des}}}
\newcommand{\QdemandIJPeak}{\dot{Q}_{\mathrm{des},ij,t_\mathrm{peak}}}
\newcommand{\QdemandPeak}{\dot{Q}_{\mathrm{des},t_\mathrm{peak}}}
\newcommand{\QdemandPeakIJ}{\dot{Q}_{\mathrm{des},t_\mathrm{peak},ij}}
\newcommand{\Qdemandi}[1]{\dot{Q}_{\mathrm{des},#1}}
\newcommand{\COPsub}{\mathrm{COP}_{\con}}	% cop hp substation
\newcommand{\COPsubIJ}{\mathrm{COP}_{\con,ij}}	% cop hp substation
\newcommand{\COPsubi}[1]{\mathrm{COP}_{\con,#1}}
\newcommand{\bigM}{M}	% big M parameter for constraints with binary variables
% producer side
\newcommand{\COPprod}{\mathrm{COP}_{\pro}}	% cop HP prod
\newcommand{\COPprodIJ}{\mathrm{COP}_{\pro,ij}}
\newcommand{\COPprodi}[1]{\mathrm{COP}_{\pro,#1}}

%% physical properties
% Hydraulic properties
\newcommand{\Reynolds}{Re}      % Reynoldsnumber
\newcommand{\density}{\rho}     % water density
\newcommand{\densityValue}{983}     % water density \kilogram\per\meter^3
\newcommand{\viscosity}{\mu}    % viscosity
\newcommand{\spHeatCap}{c_{\mathrm{p}}}
\newcommand{\spHeatCapValue}{4185} % \joule\per{\kilogram\kelvin}

\newcommand{\maxPressure}{\Delta p_\mathrm{max}}
\newcommand{\maxPressureValue}{15} % bar

%% cost terms
% J terms
\newcommand{\costFull}{\mathcal{J}}
\newcommand{\costi}[1]{\cost_{\mathrm{#1}}}
% C cost factors
\newcommand{\npvCost}{C}				 % general cost paramter capital C
\newcommand{\prodCostSpecificInv}{\npvCost_{\textrm{C,prod}}} % producer specific capex cost
\newcommand{\prodCostSpecificOper}{\npvCost_{\textrm{O,prod}}} % producer operational cost price
\newcommand{\prodCostSpecificOperTime}{\npvCost_{\textrm{O,prod,\timeVar}}} % producer operational cost price
\newcommand{\prodCostSpecificOperDT}{\npvCost_{\textrm{O,prod},\dayTimeVar}} % producer operational cost price
\newcommand{\conCostSpecificInv}{\npvCost_{\textrm{C,con}}} % consumer specific capex cost
\newcommand{\conCostSpecificInvIJ}{\npvCost_{\textrm{C,con,ij}}}
\newcommand{\conCostSpecificOper}{\npvCost_{\textrm{O,con}}} % consumer operational cost price
\newcommand{\storCostInv}{\npvCost_{\textrm{C},\stor}} % consumer specific capex cost

\newcommand{\subCAPEX}{CAP} % capex cost index
\newcommand{\subOPEX}{OP}   % opex cost index
\newcommand{\subPump}{pump}
\newcommand{\efficiency}{\eta}
\newcommand{\pumpEff}{\efficiency_{\mathrm{pump}}}
\newcommand{\CHPturbineEff}{\efficiency_{\pro,\CHPIndex}}
\newcommand{\cPumpOPEX}{\npvCost_{\mathrm{pO}}}
\newcommand{\cPumpOPEXi}{\npvCost_{\mathrm{O,pump}}}

% strings to describe the cost terms
\newcommand{\Jpump}{\text{Pump OPEX}}
\newcommand{\JheatCentral}{\text{Central heat pump OPEX}}
\newcommand{\JheatSubstation}{\text{HP substation OPEX}}
\newcommand{\JCAPEXSubstation}{\text{HP substation CAPEX}}
\newcommand{\Jtot}{\text{Total cost}}

%% environmental parameters
\newcommand{\TOutside}{\temperature_\infty}
\newcommand{\TOutsideTime}{\temperature_{\infty,\timeVar}}
\newcommand{\TGround}{\temperature_{\mathrm{ground}}}
\newcommand{\TGroundTime}{\temperature_{\mathrm{ground},\timeVar}}

\newcommand{\measNonDiscrete}{M_{nd}}

% parameters related to the pipe and its thermal model
\newcommand{\length}{L} % pipe length
\newcommand{\thermR}{U} % thermal resistance
\newcommand{\diameter}{d}
\newcommand{\subScriptOuterD}{o}  % subscript outer pipe diameter
\newcommand{\depthPipe}{h} % depth pipe
\newcommand{\depthPipeValue}{1} % depth pipe value
\newcommand{\ratioInsul}{r} % insualtion ratio
\newcommand{\ratioInsulValue}{1.87} % insualtion ratio
\newcommand{\condInsul}{\lambda_{\mathrm{i}}} % conductivity thermal insulation
\newcommand{\condInsulValue}{0.0225} % conductivity thermal insulation
\newcommand{\condGround}{\lambda_{\mathrm{g}}} % conductivity ground
\newcommand{\condGroundValue}{1} % conductivity ground
\newcommand{\thermRCorrFactor}{F_{\thermR,\mathrm{corr}}} 
\newcommand{\thermRCorrFactorValue}{0.25} 

\newcolumntype{L}[1]{>{\raggedright\arraybackslash}p{#1\columnwidth}} %MS to allow line breaks AND specific sizing in tabularx tables

\section{An automated storage sizing optimization framework for DHNs}
\label{sec:method}

In this section, the storage sizing for DHNs is formulated as a mathematical optimization problem and the solution procedure is described. The DHN simulation and optimization is based on previous developments by the authors of this paper, in particular \cite{wack2024multi}. Therefore, the focus of this section is on the new developments to enable cost-optimal storage sizing for DHNs. For completeness, the main structure of the underlying optimization problem is outlined.

\subsection{Notation}
\label{subsec:notation}

DHNs are a network technology and can be depicted with a directed graph $\dirGraph(\setNodes,\setEdges)$. Here, $\setNodes$ denotes the set of nodes, while $\setEdges$ represents the set of edges present in the graph. An overview of the graph notation used to formulate optimization problems in this paper is given in table \ref{tab:DHNnotation}.

\begin{table}[t]
	\centering
	\caption{Graph and index notation for the DHN.}
	\label{tab:DHNnotation}
	\begin{tabularx}{\columnwidth}{XX}
		\toprule
		\textbf{Notation} & \textbf{Description} \\
		\midrule
		$\dirGraph = (\setNodes, \setEdges)$ & DHN as a directed graph.\\
		$\setNodes = \Npro \cup \Ncon \cup \Nstor \cup \Njun$ & Set of nodes, including heat
		producers, consumers, storages, and junctions.\\
		$\setEdges = \Epro \cup \Econ \cup \Estor \cup \Epipe$ & Set of edges, including heat
		producers, consumers, storages and pipes.\\
		$\Econ = \Econhx \cup \Econbyp$ & Set of consumer edges, including heat
		exchanger substations and bypasses.\\
		$\Epro = \EproCHP$ & Set of producer edges, containing the central combined heat and power producer.\\
		$\Estor = \EstorHot \cup \EstorCold$ & Set of storage edges, containing the hot (supply) and cold (return) edges.\\
		$(\gi,\gj)$ or $\gi\gj$ & Directed edge going from node $\gi$ to node $\gj$. \\
		\bottomrule
	\end{tabularx}
\end{table}

\subsection{Storage sizing as an optimization problem}
\label{subsec:optimization problem}

In this work, the heat producer is assumed to be given, as are the DHN pipes. Hence, there are no producer capacity or pipe diameter design variables. For the optimal storage sizing, two new design variables are introduced for each storage. The first represents its size (time-independent), while the second represents its (dis-)charging operation (time-dependent). Additionally, a bypass is considered in this work at every consumer substation.

Designing new storages for a DHN requires to consider the performance throughout future operation. This is especially important if the potential economic benefit of the storage stems from taking advantage of varying short-term heat production prices, as in this study. Therefore, optimizing the storage sizing means solving a multi-period optimization problem. To account for the time-dependent nature of this optimization problem, a set of days $\dayVarSet=\{1,\dots,\ndays\}$ and time segments per day $\timeVarSetd=\{1,\dots,\ntsd\}$ are defined, where $\ndays$ is the number of days and $\ntsd$ is the number of time segments per day. Hence, an individual time segment is defined by the indices $\dayVar \in \dayVarSet$ and $\timeSegVar \in \timeVarSetd$. 

\begin{figure}[t]
	%singlecolumn figure
	\includegraphics[width=1\columnwidth]{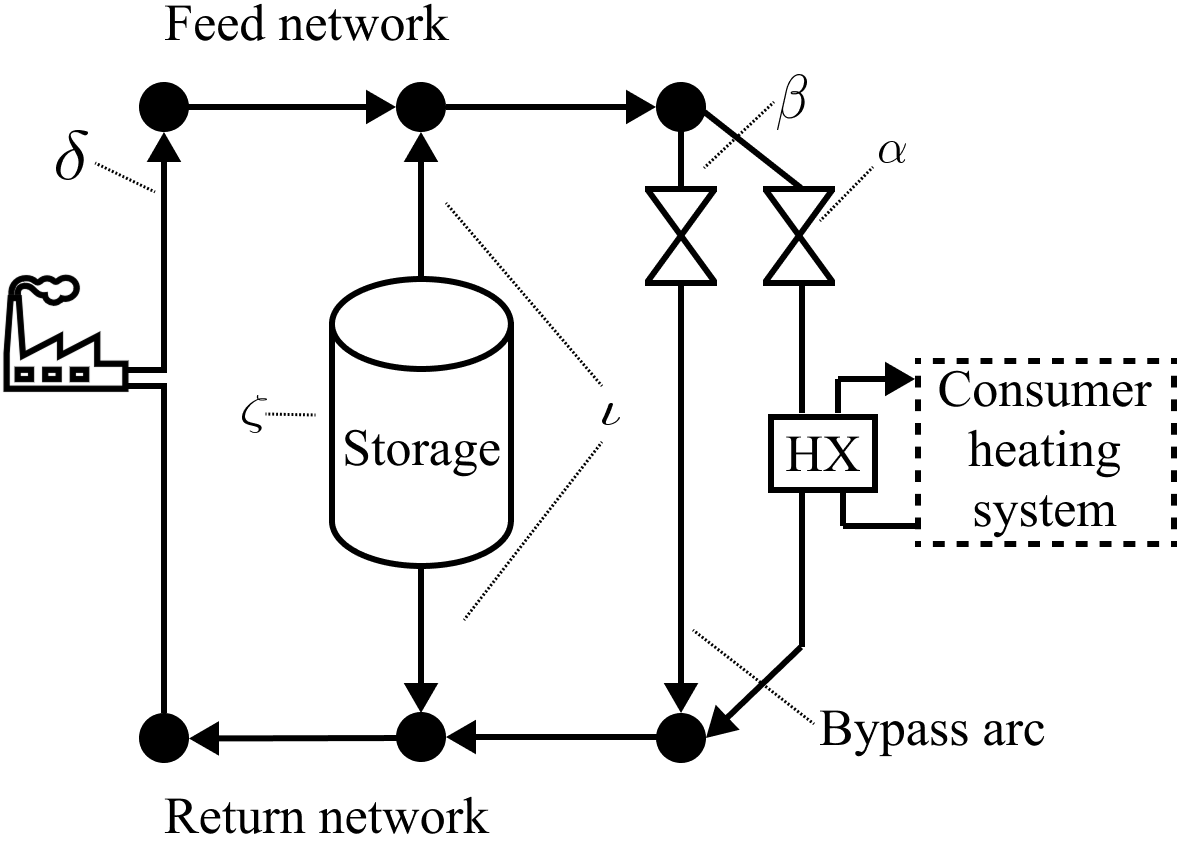}	
	\caption{District heating network with its heat producer and (potential) heat storage. The storage size is optimized together with its (dis-)charging flow and the operation of the network.}
	\label{fig:DHN_withStorage}
\end{figure}

The total discounted project cost $\costFull\left(\designVar,\stateVar\right)$ is minimized by setting for each storage its size $\ve{\storSize}$ and for all periods its (dis-)charging operation $\ve{\storFlow}$, and for the producer its mass flow rate $\ve{\prodInput}$ for all periods. To ensure a feasible operation, the operating variables (valves) of the consumer heat exchanger (HX) substations $\ve{\HXValve}$ and bypasses $\ve{\bypValve}$ must be chosen for all consumers and all periods. The role of each design variable is visualized in figure \ref{fig:DHN_withStorage}. The storage sizing optimization problem for an existing DHN is thus formulated as
\begin{equation}
	\begin{aligned} \label{eq:optProblem}
		\min_{\designVar,\stateVar} &\qquad
		\costFull \left(\designVar,\stateVar\right) &\\
		s.t.& \qquad \defModelConstraintsDT = 0,\quad\forall \dayVar \in \dayVarSet\,,\forall \timeSegVar \in \timeVarSetd\,, &\\
		& \qquad \defStateConstraintsDT \leq 0,\quad\forall \dayVar \in \dayVarSet\,,\forall \timeSegVar \in \timeVarSetd, & \\
		& \qquad \opVarBoxConstraints{\designVar}\,, &
	\end{aligned}
\end{equation}
where the design variables $\defDesignVar$ and the physical variables $\defStateVar$ are optimized to satisfy the model equations $\defModelConstraintsDT = 0$. Note that all variables are defined as continuous variables, making the optimization problem a continuous, nonlinear problem (NLP). For ease of notation, the time-dependent variable vectors $\designVarDT$ and $\stateVarDT$ of all periods are combined in the vectors $\designVar$ and $\stateVar$, respectively.

\begin{figure}[t]
	%singlecolumn figure
	\includegraphics[width=1\columnwidth]{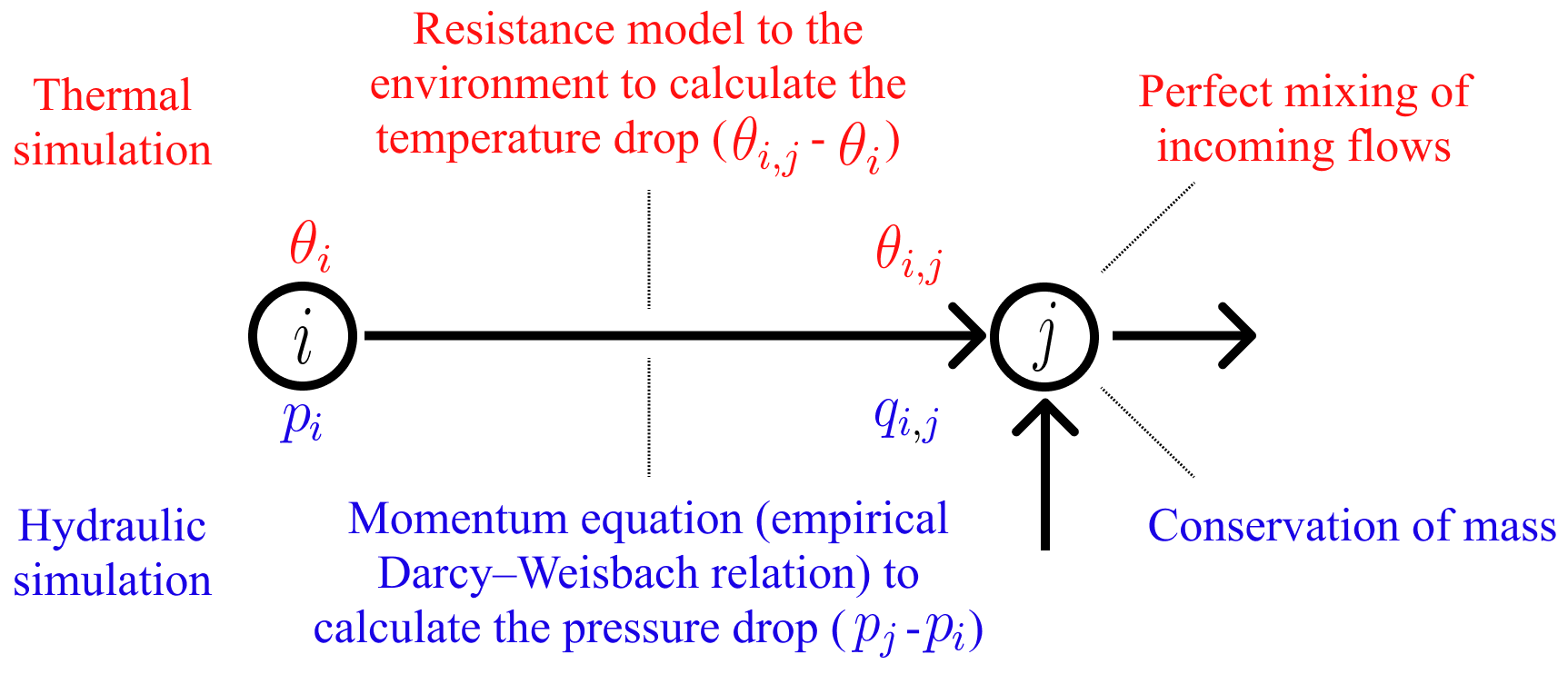}	
	\caption{Thermal and hydraulic simulation of the DHN, accounting for heat and
		momentum losses.}
	\label{fig:simulationGraph}
\end{figure}

The set of nonlinear model equations accurately represents the hydraulic and thermal transport problem in the network and heat storage, accounting for momentum and heat losses, to model pressure and temperature drops throughout the network. It includes the flow rates $\defFlow$, nodal pressures $\defPressure$, and nodal and pipe exit temperatures $\defTemp$. The temperature vector $\defTemp = \ve{\temperature} - \TOutside$ is defined as the difference between the absolute water temperature vector $\ve{\temperature}$ and the outside air temperature $\TOutside$. Figure \Ref{fig:simulationGraph} shows the core parts of the hydraulic and thermal simulation of the DHN. For the remainder of the paper the notations $i,j$ and $ij$ are used interchangeably to refer to a network edge. This work introduces a new state variable, the storage mass $\ve{\stateVarStorMass}$, which will be explained further in the following sections. For more details on the DHN simulation, the interested reader is referred to \cite{wack2024multi}.

Additional technological constraints are represented by state constraints $\defStateConstraintsDT$. Among them is a constraint on the consumer heat demand satisfaction, this is formulated as 
\begin{multline}\label{eq:ThermalComfort}
	\frac{\Qdemandi{\gi\gj,\dayTimeVar}-\stateVarHeat_{\gi\gj,\dayTimeVar}}{\Qdemandi{\gi\gj,\dayTimeVar}} = 0 \quad \forall ij \in \Econhx ,
\end{multline}
where $\Qdemandi{\gi\gj,\dayTimeVar}$ is the heat demand and $\stateVarHeat_{\gi\gj,\dayTimeVar}$ the transferred heat. This equality constraint is imposed via a lower and an upper bound inequality constraint.

Another constraint is imposed on the maximal pressure increase, defined on each producer edge which can be stated as
\begin{multline}\label{eq:MaxPressure}
	\frac{(\stateVarPress_{j,\dayTimeVar}-\stateVarPress_{i,\dayTimeVar})-\maxPressure}{\maxPressure} \, \leq 0 \quad \forall ij \in \Epro .
\end{multline}
The other technological constraints are newly introduced in this work and are described in the next sections.

\newcommand{\npv}{NPV}
\newcommand{\npvConst}{f}
\newcommand{\npvN}{A}
\newcommand{\npvNValue}{20} % considered life time (years)
\newcommand{\npvt}{l}
\newcommand{\npvDiscount}{e}
\newcommand{\npvDiscountValue}{5} % discount rate (%)
\newcommand{\npvConstOPEX}{\npvConst_{\mathrm{\subOPEX}}}
\newcommand{\periodWeight}{\omega}

The discounted lifetime cost of the storage sizing for a DHN $\costFull$ includes both CAPEX and OPEX related to the storage, the central heat producer and the circulation pumps and is defined as: 
\begin{multline}\label{eq:totalCost}
	\costFull\left(\designVar,\stateVar\right) =  \sum_{\gi\gj \in \Estor}\costi{\stor,\subCAPEX,\textit{ij}} \\
	+ \npvConstOPEX \sum_{\dayVar = 1}^{\ndays} \sum_{\timeSegVar = 1}^{\ntsd} \periodWeight_{\dayTimeVar} \Bigg[ \sum_{\gi\gj \in \Estor} \costi{\subPump,\subOPEX, \textit{ij},\textit{\dayTimeVar}} \\
	+ \sum_{\gi\gj \in \Epro} [\costi{\pro,\subOPEX}
	+\costi{\subPump,\subOPEX}]_{\textit{ij},\textit{\dayTimeVar}} \Bigg]\,,
\end{multline}
with $\npvConstOPEX = \sum_{\npvt=1}^{\npvN} \frac{1}{\left(1+\npvDiscount\right)^\npvt} $ assuming an investment horizon of $\npvN = \npvNValue\,\mathrm{years}$ and a discount rate of $\npvDiscount=\npvDiscountValue\%$. Here, it is assumed that the chosen days $\dayVar$ and time segments $\timeSegVar$ are sufficiently representative for a full year operation and $\sum_{\dayVar = 1}^{\ndays} \sum_{\timeSegVar = 1}^{\ntsd} \periodWeight_{\dayTimeVar} = 1$. 

The operational cost of the distribution pumps at the heat producers is defined as in \cite{SOLLICH2025Producer} and computed with
\begin{multline}\label{eq:PumpopexProducer}
	\costi{\subPump,\subOPEX, \textit{ij}, \textit{\dayTimeVar}} = \frac{\KOPEX}{\pumpEff} \cPumpOPEXi \, \stateVarFlow_{ij,\dayTimeVar} \, (\stateVarPress_{j,\dayTimeVar}-\\ \stateVarPress_{i,\dayTimeVar}) \quad \forall ij \in \Epro ,
\end{multline}
while $\cPumpOPEXi$ is the electricity purchase price and the conversion factor $\KOPEX$ defines the network's number of operating hours per year. The other individual cost components $\costi{}$ are newly introduced and described in the next sections.

\subsection{Storage time-dependency and time series aggregation}
\label{subsec:time-dependency and time series aggregation}

To ensure computational tractability of the proposed optimization method, temporal aggregation of the time series heat demand, outside temperature, and electricity wholesale price is used. More specifically, an aggregation into representative days (of 24 hours length) and variable-length time segments per day is used, similar to \cite{Resimont2021, CERUTI2025}. For the clustering, the tsam Python package \cite{tsamClustering_webpage} was used, introduced in \cite{KOTZUR2018_tsamIntro}. 

\begin{figure}[t]
	\centering
	%singlecolumn figure
	\includegraphics[width=\columnwidth]{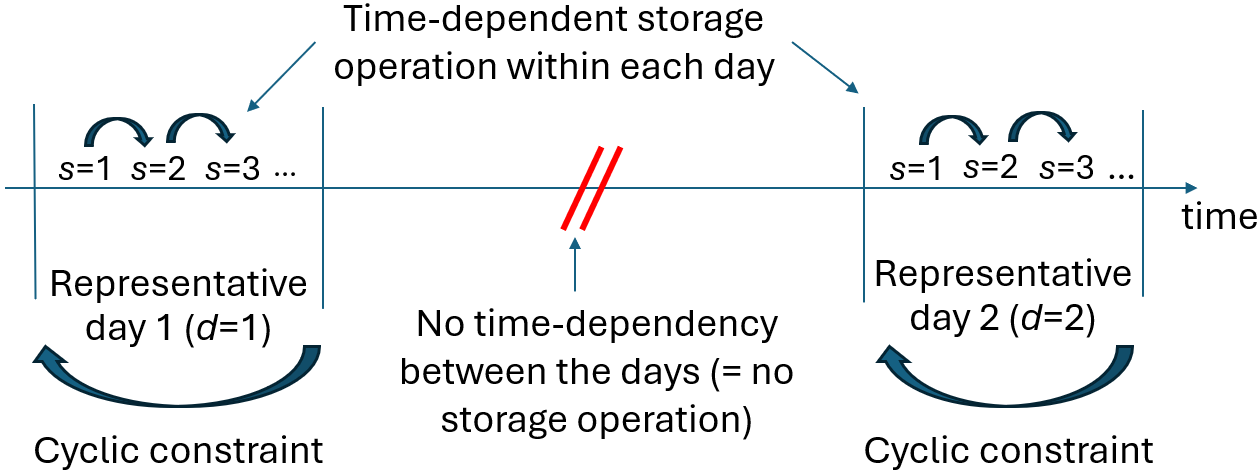}	
	\caption{The operational regime of the short-term storage considered in this work and its link to the time-dependency.}
	\label{fig:timeDependencyApproach}
\end{figure}

Since this work considers only short-term heat storage, the (dis-)charging of the storage is only permitted within a representative day, while a cyclic constraint ensures that the energy content at the end of a day is equal to that at the beginning of the day, see figure \ref{fig:timeDependencyApproach}. Hence, the representative days can be simulated in parallel as no time-dependency exists between them.   

The DHN is modeled with a quasi-steady state approach, as it was proposed in \cite{wack2024multi}. Hence, at each time segment, a steady state in the DHN is assumed and the hydraulic and thermal state is assumed to be independent of the states from the other time segments. This can be justified since the considered time segments in this work represent multiple hours. In comparison to this, hydraulic transients in the DHN can be considered as instantaneous \cite{STEVANOVIC2007}. Thermal transients in DHN can last in principle for longer, even up to several hours \cite{BENONYSSON1995}. However, since the proposed methodology assumes a steady-state at the consumer heating system and substation at each time segment and since the producer supply temperature is constant in the considered case, the quasi-steady state approach can also be justified for the DHN's thermal state.

This hybrid approach, which models thermal transients only where necessary (i.e., in the storage), is one of the core pillars of creating a scalable, physics-based (nonlinear) design optimization framework for storage in DHNs. 

\subsection{DHN storage integration}
\label{subsec:storageMethodSec}

\begin{figure}[t]
	\centering
	%singlecolumn figure
	\includegraphics[width=\columnwidth, height=5cm, keepaspectratio]{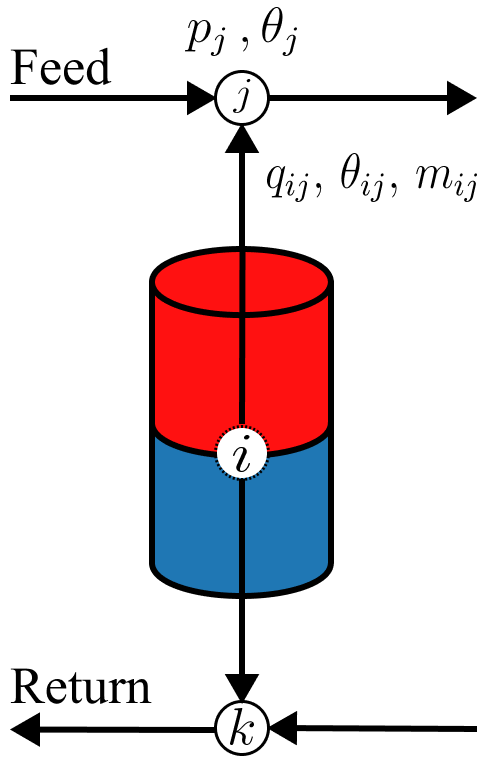}	
	\caption{The 2-zone storage model and its integration into the DHN superstructure\protect\footnotemark. Note that the shown state variables are analogically defined for the cold (blue) zone of the storage.}
	\label{fig:storage2zoneModel}
\end{figure}

In this work, a cylindrical water tank is considered for the storage, which is installed above the ground. The storage is modeled with two perfectly mixed (hot and cold) zones that vary in temperature and size (and therefore mass) over time but that cannot exchange heat with each other, which renders a 1D model, see figure \ref{fig:storage2zoneModel}. Hence, the model assumes a perfectly stratified tank with a perfect thermocline between the two zones. Such 1D models are commonly used when modeling a heat storage as part of a larger energy system as they combine reasonable computational complexity with sufficient accuracy under transient conditions \cite{JODEIRI2024}. A 2-zone storage model was already used for example in \cite{Vandermeulen2020PhDThesis} and \cite{Delubac2023}, for optimal control and optimal design of DHNs, respectively.
\footnotetext{Note that the node $i$ inside the storage is not part of the superstructure as no variables are assigned to it, it is only shown to have the edges $ij$ and $ik$ well-defined between two nodes.}
In order to model the storage in the optimization algorithm, the optimization's superstructure (the network graph) is extended by one node and one edge for the hot and the cold zone, respectively, as shown in figure \ref{fig:storage2zoneModel}. Moreover, with the storage mass $\stateVarStorMass$, defined for each zone, a new state variable is introduced. 
 
The time-dependent mass and energy balance per zone result in differential equations which are discretized with an Euler forward scheme. This allows that all state variables can be defined for each day $\dayVar$ and each time segment $\timeSegVar$. However, for reasons of compactness and better representation of the underlying physics, the storage model equations are written in non-discretized form in the following sections, which means for example that $\stateVarFlow_{ij}(t)$ is used instead of $\stateVarFlow_{ij,\dayTimeVar}$ to represent the time-dependent flow rate of the edge $ij$.

\subsubsection{Storage model - hydraulic equations}
\label{subsubsec:storageModel_hydrEquations}

The storage (dis-)charging operation is optimized through the operational variable $\storFlow_{ij}(\timeVar)$ which represents a flow rate. This operational variable is linked to the flow rate $\stateVarFlow_{ij}(\timeVar)$ through a boundary condition which reads
\begin{equation}\label{eq:storHydrFlowBC_I}
	\stateVarFlow_{ij}(\timeVar) = \storFlow_{ij}(\timeVar) \quad \forall ij \in \EstorHot,
\end{equation}
while negative values represent an inflow into the hot zone, positive values an outflow. The same flow rate value is also imposed on the cold zone but with the opposite sign:
\begin{multline}\label{eq:storHydrFlowBC_II}
	\stateVarFlow_{ij}(\timeVar) = -\stateVarFlow_{ik}(\timeVar) \quad \forall ij \in \EstorHot,\, ik \in \EstorCold.
\end{multline}
Thereby, it is ensured that the total mass of the storage is time-independent. 

The hot zone storage mass $\stateVarStorMass_{ij}(t)$ can be determined from mass conservation on the hot zone which reads
\begin{equation}\label{eq:storHydrMass_diff}
	\frac{d\stateVarStorMass_{ij}(\timeVar)}{d\timeVar} = -\stateVarFlow_{ij}(\timeVar)\density \quad \forall ij \in \EstorHot,
\end{equation}
with $\density$ being the density of the water. The initial hot mass at the beginning of each day is imposed as a fixed boundary condition.

The storage mass in the cold zone $\stateVarStorMass_{ik}(t)$ can be determined from the total storage mass $\stateVarStorMass_{\totalMassIDX,ij}$ and the hot zone storage mass which reads
\begin{multline}\label{eq:storHydrMass_cold}
	\stateVarStorMass_{ik}(t) = \stateVarStorMass_{\totalMassIDX,ij} - \stateVarStorMass_{ij}(t) \quad \forall ik \in \EstorCold, \\ \, ij \in \EstorHot.
\end{multline}

The total storage mass can be calculated from the storage volume which is given by the normalized storage size design variable $\storSize_{ij}$ multiplied by the maximal storage volume $\storMaxVolume$: 
\begin{multline}\label{eq:storHydrMass_total}
	\stateVarStorMass_{\totalMassIDX,ij} = \storSize_{ij}\storMaxVolume\density \quad \forall ij \in \EstorHot.
\end{multline}

\subsubsection{Storage model - thermal equations}
\label{subsubsec:storageModel_thermEquations}

\newcommand{\QstorLossIJ}{\dot{Q}_{\mathrm{loss},ij}}
\newcommand{\QstorChargingIJ}{\dot{Q}_{\mathrm{char},ij}}
\newcommand{\QstorDischargingIJ}{\dot{Q}_{\mathrm{dischar},ij}}

To determine the hot zone storage temperature $\stateVarTemp_{ij}(\timeVar)$, an energy balance is imposed on the hot zone which reads
\begin{multline}\label{eq:storTherm_HotZone_diff}
	\spHeatCap\frac{d\left(\stateVarStorMass(\timeVar)\stateVarTemp(\timeVar)\right)_{ij}}{d\timeVar} =
	-\QstorLossIJ(\timeVar) \\
	+ \QstorChargingIJ(\timeVar)\,
	- \QstorDischargingIJ(\timeVar)	
	\quad \forall ij \in \EstorHot,
\end{multline}
with $\spHeatCap$ being the specific heat capacity of water, $\QstorLossIJ(\timeVar)$ the storage heat loss, $\QstorChargingIJ(\timeVar)$ the charged energy, and $\QstorDischargingIJ(\timeVar)$ the discharged energy. 

The hot zone heat loss to the surrounding (ambient) air is calculated via an overall heat transfer coefficient $\Ustor$ and the surface area of the hot zone $\storSurface_{ij}(t)$ which reads
\begin{multline}\label{eq:storTherm_Qloss_diff}
	\QstorLossIJ(\timeVar) = \Ustor\,\storSurface_{ij}(\timeVar)\stateVarTemp_{ij}(\timeVar)
	\quad \forall ij \in \EstorHot.
\end{multline}
The hot zone surface area can be derived from the hot zone storage mass, the total storage volume, and its height-diameter ratio\footnote{The height-diameter ratio of the storage tank is assumed to be 3.}. The energy charged into the hot zone can be calculated from
\begin{multline}\label{eq:storTherm_Qchar_diff}
	\QstorChargingIJ(\timeVar) = -\min(\stateVarFlow_{ij}(\timeVar),0)\density\spHeatCap\stateVarTemp_{j}(\timeVar)\\
	\quad \forall ij \in \EstorHot,
\end{multline}
with $\stateVarTemp_{j}(\timeVar)$ being the DHN supply temperature at the storage. The discharged energy from the hot zone can be determined from
\begin{multline}\label{eq:storTherm_Qdischar_diff}
	\QstorDischargingIJ(\timeVar) =  \max(\stateVarFlow_{ij}(\timeVar),0)\density\spHeatCap\stateVarTemp_{ij}(\timeVar)\\
	\quad \forall ij \in \EstorHot.
\end{multline}
The equations \ref{eq:storTherm_Qloss_diff} - \ref{eq:storTherm_Qdischar_diff} are analogically defined for the cold zone of the storage. The temperatures of both the hot and the cold zone at the first time segment of each day are imposed as fixed boundary conditions.

The amounts of energy charged and discharged depend on the direction of the storage flow rate $\stateVarFlow_{ij}(\timeVar)$, which is inherently a discrete decision: at any given time $\timeVar$, the storage is either charging or discharging. This directional dependency requires the use of $\min$ and $\max$ functions in equations \ref{eq:storTherm_Qchar_diff} and 
\ref{eq:storTherm_Qdischar_diff}, introducing non-differentiabilities into the optimization problem that pose a fundamental challenge for gradient-based solvers. Such non-differentiabilities are the core difficulty of problems that are intrinsically mixed-integer nonlinear 
programs (MINLPs). To overcome this, the $\min$ and $\max$ functions are replaced by the \textit{softplus} function \cite{Dugas2000_softPlus}, a smooth approximation that eliminates the non-differentiabilities and allows the problem to be formulated as a continuous NLP, solvable efficiently via gradient-based optimization.

\subsubsection{Storage constraints}
\label{subsubsec:storageConstraints}

In addition to satisfying the physical storage model (flow and heat transfer equations), constraints are defined for the storage to ensure that a useful optimization problem is solved. First, since the 2-zone storage model is only well-defined if the shares of the hot and cold zone mass remain between 0 and 1, the hot zone mass share is constrained which reads
\begin{equation}\label{eq:storHydrMass_constraint}
	0.05 \leq \frac{\stateVarStorMass_{ij}(\timeVar)}{\stateVarStorMass_{\totalMassIDX,ij}} \leq 0.95 \quad \forall ij \in \EstorHot.
\end{equation}
Here, 5 and 95 \% are used as bounds to ensure a robust simulation. For the same reason, this constraint is imposed through an interior-point method and not trough the Augmented-Lagrangian to ensure that the mass shares remain within these bounds at every iteration of the optimization.
Through equations \ref{eq:storHydrMass_cold} and \ref{eq:storHydrMass_constraint}, the cold zone mass share is constrained equivalently. 

Second, as mentioned in section \ref{subsec:time-dependency and time series aggregation}, a cyclic constraint is imposed to ensure that the energy at the end of each representative day equals its initial energy which reads
\begin{multline}\label{eq:storCyclic_constraint}
	\left(\stateVarStorMass_{\initINDEX}\stateVarTemp_{\initINDEX}\right)_{ij,\dayVar} + \left(\stateVarStorMass_{\initINDEX}\stateVarTemp_{\initINDEX}\right)_{ik,\dayVar} = \\
	\left(\stateVarStorMass\stateVarTemp\right)_{ij,\dayVar,\timeSegVar=\ntsd+1} + \left(\stateVarStorMass\stateVarTemp\right)_{ik,\dayVar,\timeSegVar=\ntsd+1} \\
	\quad \forall ij \in \EstorHot,\, ik \in \EstorCold.
\end{multline}
note that $\stateVarStorMass$ and $\stateVarTemp$ are not solved in the simulation for $\timeSegVar=\ntsd+1$ but that they are implicitly determined by the discretized energy balance from $\timeSegVar=\ntsd$ to $\timeSegVar=\ntsd+1$, based on equation \ref{eq:storTherm_HotZone_diff}.

Lastly, the maximal pressure increase of the circulation pump at the storage is constrained, same as for the heat producer circulation pump, which can be stated as
\begin{multline}\label{eq:MaxPressureStorage}
	\frac{(\stateVarPress_{j,\dayTimeVar}-\stateVarPress_{l,\dayTimeVar})-\maxPressure}{\maxPressure} \, \leq 0 \\ \quad \forall ij \in \EstorHot,\, ik \in \EstorCold.
\end{multline}

\subsubsection{Storage costs}
\label{subsubsec:storageCosts}

\begin{figure}[t]
	\centering
	%singlecolumn figure
	\includegraphics[width=\columnwidth]{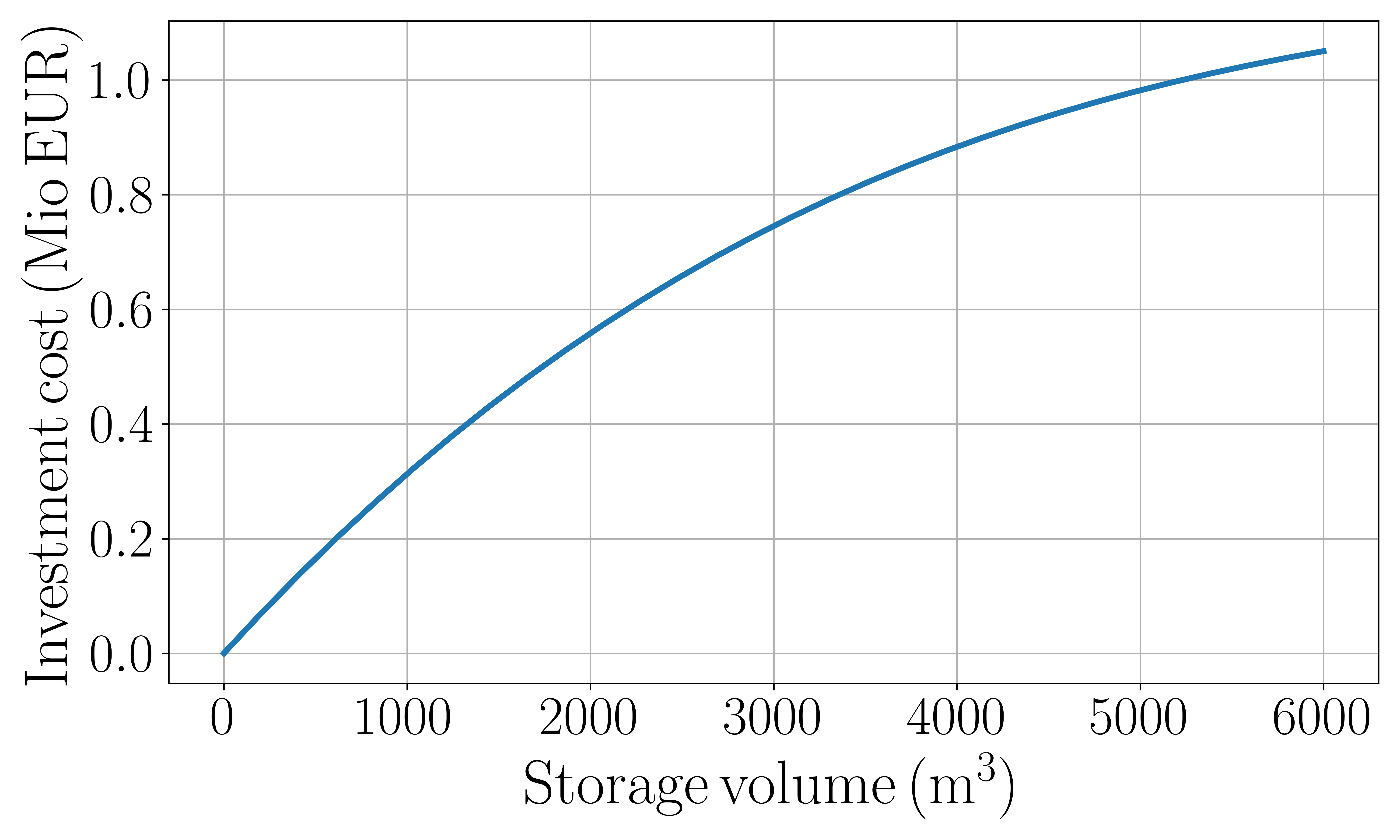}	
	\caption{Fitted investment cost $\storCostInv$ of the heat storage.}
	\label{fig:storage_CAPEX_plot}
\end{figure}

The OPEX cost of the circulation pump at the storage is considered analogically to the OPEX of the producer circulation pump and can be computed with the (dis-)charging flow rate $\stateVarFlow_{ij,\dayTimeVar}$ and the pressure difference between the supply and return side of the DHN at the storage $(\stateVarPress_{j,\dayTimeVar}-\stateVarPress_{l,\dayTimeVar})$ which reads 
\begin{multline}\label{eq:PumpopexStorage}
	\costi{\subPump,\subOPEX, \textit{ij}, \textit{\dayTimeVar}} = \frac{\KOPEX}{\pumpEff} \cPumpOPEXi \, \lvert\stateVarFlow_{ij,\dayTimeVar}\rvert \, (\stateVarPress_{j,\dayTimeVar}-\\ \stateVarPress_{l,\dayTimeVar}) \quad \forall ij \in \EstorHot,\, ik \in \EstorCold.
\end{multline}
For the storage investment cost (CAPEX), a fit of the total investment cost as a function of the storage size is derived, $\storCostInv(\storSize_{ij} \, \storMaxVolume)$, which is visualized in figure \ref{fig:storage_CAPEX_plot}, showing a decrease in the specific investment cost with increasing storage size. The CAPEX can then be calculated via
\begin{multline}\label{eq:storage_CAPEX}
	\costi{\stor,\subCAPEX,\textit{ij}} =  \storCostInv(\storSize_{ij} \, \storMaxVolume) \\ \quad \forall ij \in \EstorHot,
\end{multline}
with $\storSize_{ij}$ being the normalized storage volume design variable. More details on the fit and its underlying data are provided in table \ref{tab:storageCAPEXfit}.

\subsection{Combined heat and power OPEX}
\label{subsec:CHP_OPEX}

The CHP unit at the studied DHN is an extraction CHP unit, where steam is extracted for the DHN heat supply before it enters the turbine, thereby reducing the amount of electricity produced by the turbine. Hence, the operational cost (OPEX) of heat provision can be defined as the loss in electricity sales revenue when extracting heat for the DHN which reads
\begin{multline}\label{eq:chpProdOPEX}
	\costi{\pro,\subOPEX,\textit{ij},\textit{\dayTimeVar}} = \KOPEX\,\prodCostSpecificOperDT\,\CHPturbineEff\,\density\,\spHeatCap\left(\stateVarFlow\,\Delta\theta\right)_{ij,\dayTimeVar} \\ \quad \forall ij \in \EproCHP.
\end{multline}
Here, the electricity wholesale price is given by $\prodCostSpecificOperDT$, while $\density\,\spHeatCap(\stateVarFlow\,\Delta\theta)_{ij,\dayTimeVar}$ represents the heat supplied to the DHN, and $\CHPturbineEff$ is the efficiency of the turbine\footnote{The presented CHP OPEX is not the actual (contractual) heat cost of the CHP plant, which cannot be disclosed publicly. It represents a generic OPEX term for extraction CHP units.}.

\subsection{Optimization methodology}
\label{subsec:optimization methodology}

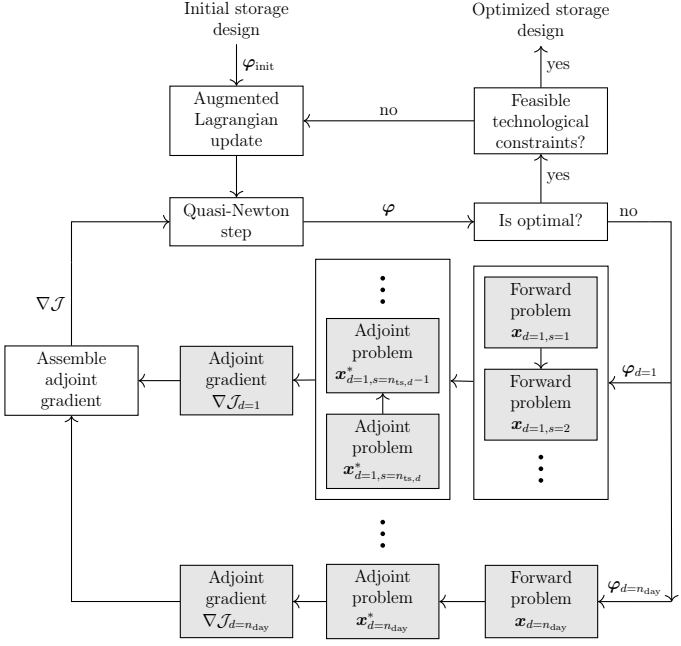
\begin{figure}[t]          
	\centering
	%singlecolumn figure  
	\begin{tikzpicture}[scale=0.55, every node/.style={transform shape},
		node distance=1.7cm]
		
		\tikzstyle{block} = [rectangle, draw, minimum width=3.2cm, minimum height=0.9cm, align=center]
		\tikzstyle{block_noFrame} = [rectangle, minimum width=3.2cm, minimum height=0.9cm, align=center]
		\tikzstyle{smallblock} = [rectangle, draw, minimum width=2.7cm, minimum height=0.9cm,
		align=center, fill=gray!20]
		\tikzstyle{arrow} = [->, thin]
		\tikzstyle{line} = [-, thin]
		\tikzstyle{lineDashed} = [dashed, thin]
		
		% Top row
		\node[block] (aug0) {Augmented\\Lagrangian\\update};
		\node[block_noFrame, above=1cm of aug0] (init) {Initial storage\\design};
		\node[block, right=4.1cm of aug0] (feasible) {Feasible\\ technological\\constraints?};
		\node[block_noFrame, above=1cm of feasible] (optNet) {Optimized storage\\design};
		
		\draw[arrow] (init) -- node[right] {$\designVar_{\mathrm{init}}$} (aug0);
		\draw[arrow] (feasible) -- node[right] {yes} (optNet);
		\draw[arrow] (feasible) -- node[above] {no} (aug0);
		
		% Quasi-Newton block
		\node[block, below=1cm of aug0] (quasi) {Quasi-Newton\\step};
		\draw[arrow] (aug0) -- (quasi);
		
		% Optimal? block
		\node[block, right=4.1cm of quasi] (optimal) {Is optimal?};
		\draw[arrow] (quasi) -- node[above] {$\designVar$} (optimal);
		\draw[arrow] (optimal) -- node[right] {yes} (feasible);
		
		%% Forward/adjoint chain (first layer)
		% gradient nodes
		\node[smallblock, below=2.4cm of quasi, xshift=0cm] (grad1) 	{Adjoint\\gradient\\$\nabla\costFull_{\dayVar=1}$};
		% adjoint nodes
		\node[smallblock, right=0.8cm of grad1, yshift=0.6cm] (adjD1Tn1) 	{Adjoint\\problem\\$\stateVar^*_{\dayVar=1,\timeSegVar=\ntsd-1}$};
		\node[smallblock, below=0.5cm of adjD1Tn1] (adjD1Tn) 	{Adjoint\\problem\\$\stateVar^*_{\dayVar=1,\timeSegVar=\ntsd}$};
		\node[above=0.2cm of adjD1Tn1] (dotsAD1) {
			\tikz{
				\fill (0,0) circle (1.6pt);
				\fill (0,-0.28) circle (1.6pt);
				\fill (0,-0.56) circle (1.6pt);
			}
		};
		
	    % frame around the adjoint solve of day 1
		\node[draw, rectangle,
		fit = (dotsAD1) (adjD1Tn1) (adjD1Tn),
		inner sep=4pt,
		] (FrameADd1) {};
		% arrows withing the adjoint solve of day 1
		\draw[arrow] (adjD1Tn) -- (adjD1Tn1);
		
		% forward day 1
		\node[smallblock, below=0.02cm of adjD1Tn1,xshift=3.8cm,yshift=0.6cm] (fwdD1T2) 	{Forward\\problem\\$\stateVar_{\dayVar=1,\timeSegVar=2}$};
		\node[smallblock, above=0.5cm of fwdD1T2] (fwdD1T1) {Forward\\problem\\$\stateVar_{\dayVar=1,\timeSegVar=1}$};
		\node[below=0.2cm of fwdD1T2] (dotsFD1) {
			\tikz{
				\fill (0,0) circle (1.6pt);
				\fill (0,-0.28) circle (1.6pt);
				\fill (0,-0.56) circle (1.6pt);
			}
		};
		
		% frame around the forward solve of day 1
		\node[draw, rectangle,
		fit = (fwdD1T1) (fwdD1T2) (dotsFD1),
		inner sep=4pt,
		] (FrameFWDd1) {};
		% arrows withing the forward solve of day 1
		\draw[arrow] (fwdD1T1) -- (fwdD1T2);
		
		%helper coordinate to the right of is optimal
		\coordinate (optimalRightHelper) at ($(optimal.east)+(1cm,0)$);
		%helper coordinate to the right of forward d1
		\coordinate (fd1Helper) at ($(FrameFWDd1.east)+(1.5cm,0)$);
		
		\draw[line] (optimal) -- node[above]{no} (optimalRightHelper);
		\draw[line] (optimalRightHelper) -| (fd1Helper);
		\draw[arrow] (fd1Helper) -- node[above]{$\designVar_{\dayVar=1}$} (FrameFWDd1);
		\draw[arrow] (FrameFWDd1) -- (FrameADd1);
		\draw[arrow] (FrameADd1) -- (grad1);
		
		% d-th layer
		\node[smallblock, below=3.6cm of grad1] (gradD) {Adjoint\\gradient\\$\nabla\costFull_{\dayVar=\ndays}$};
		\node[smallblock, right=0.8cm of gradD] (adjD) {Adjoint\\problem\\$\stateVar^*_{\dayVar=\ndays}$};
		\node[smallblock, right=1.1cm of adjD] (fwdD) {Forward\\problem\\$\stateVar_{\dayVar=\ndays}$};
		
		%helper coordinate to the right of forward d1
		\coordinate (fwdDHelper) at ($(fwdD.east)+(1.77cm,0)$);
	
		\draw[arrow] (fwdD) -- (adjD);
		\draw[arrow] (adjD) -- (gradD);
		\draw[arrow] (fd1Helper) -- (fwdDHelper);
		\draw[arrow] (fwdDHelper) -- node[above]{$\designVar_{\dayVar=\ndays}$} (fwdD);
		
		% Assemble gradient
		\node[block, left=1.0cm of grad1] (assemble) {Assemble\\adjoint\\gradient};
	    %helper coordinate above assemble
		\coordinate (assembleHelper) at ($(assemble.north)+(0,+2cm)$);
		
		\draw[arrow] (grad1) -- (assemble);
		\draw[arrow] (gradD) -| (assemble);
		
		% Back to Quasi Newton
		\draw[line] (assemble) -- node[left] {$\nabla\costFull$} (assembleHelper);
		\draw[arrow] (assembleHelper) |- (quasi);
		
		% vertical dots above the dashed line
		\node[below=0.6cm of adjD1Tn] (dotsDaysAbove) {
			\tikz{
				\fill (0,0) circle (1.6pt);
				\fill (0,-0.28) circle (1.6pt);
				\fill (0,-0.56) circle (1.6pt);
			}
		};
				
	\end{tikzpicture}
    \caption{Flowchart of the proposed optimization framework. The algorithm efficiently decouples individual days, considers technological constraints using an Augmented Lagrangian approach, and solves the resulting NLP using a Quasi-Newton approach and adjoint gradients. The adjoint variables are represented by $\stateVar^*$. Figure adapted from \cite{wack2024multi}.}
\label{fig:tikzOptimizationFlowchart}
\end{figure}

The nonlinear optimization problem is solved using an Augmented Lagrangian approach combined with a Quasi–Newton method based on adjoint gradients. The adjoint gradients of the different representative days can be calculated in parallel as no time-dependency between the representative days exists. To obtain the day-specific gradients, the physical and adjoint states of the time segments within each day must be solved. Since a forward Euler method is used for the time discretization, the physical state (forward problem) is solved forward in time and the adjoint state (adjoint problem) backward in time. The optimization routine is visualized in figure \ref{fig:tikzOptimizationFlowchart}.

Hessian information for the Quasi–Newton method is retrieved using a BFGS algorithm. Optimality of the solution is ensured by requiring that the KKT optimality conditions are satisfied. The interested reader is referred to \citet{Nocedal2000} for more details on the mentioned algorithms. The optimization method is part of an in-house DHN design tool called PATHOPT.

\section{Optimal storage integration into a district heating network - a case study}
\label{sec:case study}
% consumer names

This section applies the proposed storage sizing methodology to a CHP-DHN in Belgium. It compares the total cost (CAPEX and OPEX) of the DHN with storage to the reference case without storage, and it analyzes the CHP and storage operations. Next, the impact of the storage size on total cost is examined in greater depth through a parameter scan. Finally, the proposed methodology is compared to a commonly used simplified approach that aggregates the entire DHN into a single demand point, highlighting the resulting differences in storage size and costs.

\subsection{The existing district heating network}
\label{subsec:CaseSetup}

\begin{figure}[t]
	\centering
	%singlecolumn figure
	\includegraphics[width=\columnwidth]{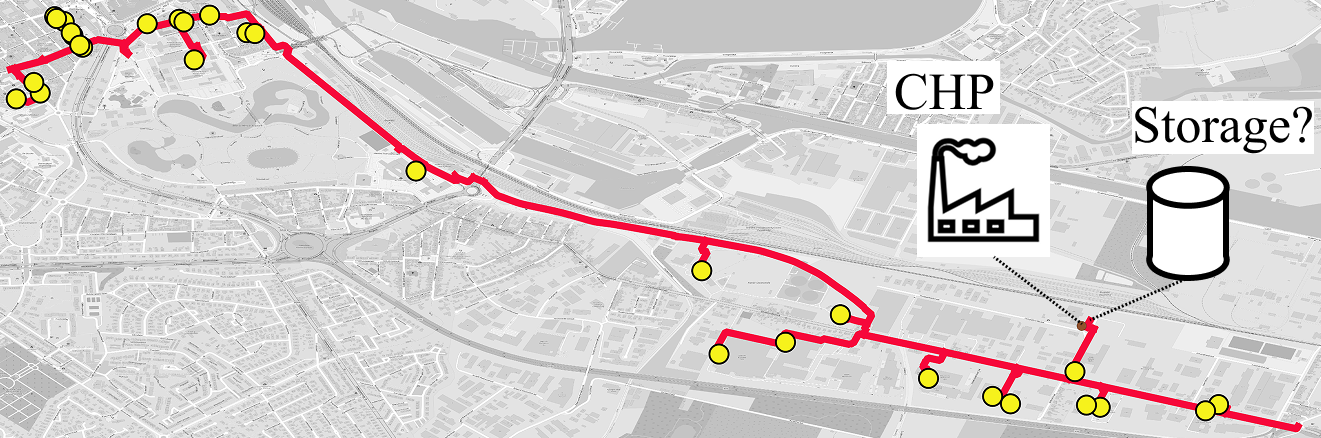}	
	\caption{The studied DHN in Oostende, Belgium. The yellow dots represent the consumers, the red lines the DHN pipes.}
	\label{fig:GISplot_caseSetup}
\end{figure}

The studied case is a 3rd generation DHN in the city of Oostende, Belgium with 30 consumers. Among the consumers are industrial clients, offices, hospitals, and apartment buildings. The total peak demand amounts to $13.4\unit{\,\mega\watt}$. The DHN has one heat source, a waste incineration CHP unit, which supplies heat at $90\unit{\,\degreeCelsius}$. There is currently no heat storage present in the network, neither centrally nor locally at the buildings. The network is visualized in figure \ref{fig:GISplot_caseSetup}.

To benefit from varying electricity wholesale prices, the DHN operator wants to investigate the value of adding a central heat storage next to the CHP unit and how the storage should be optimally sized. The investigated storage is a cylindrical water tank and is intended for short-term heat storage. The CHP unit is an extraction unit where heat for the DHN is extracted from the steam circuit before it enters the turbine, thereby reducing the electricity production and its revenue. Hence, the heat production cost could be reduced by shifting the heat extraction to moments of low electricity wholesale prices \cite{GUELPA2019}.

For the outside temperature $\TOutside$, a time series was obtained from \cite{OutsideTemperature} for the year of 2024. For the electricity wholesale price a time series of the day-ahead price for Belgium from 2024 was obtained from \cite{electricitySalePriceBeauvent}. The individual heat demand time series of all consumers were obtained from measurement data from substations meters, again from 2024. In this work, only the combined heat demand of all consumers is shown due to data privacy. However, note that the individual heat demand time series of all consumers were considered for the case study and also their individual temperature requirements.

\begin{figure}[t]
	\centering
	%singlecolumn figure
	\includegraphics[width=\columnwidth]{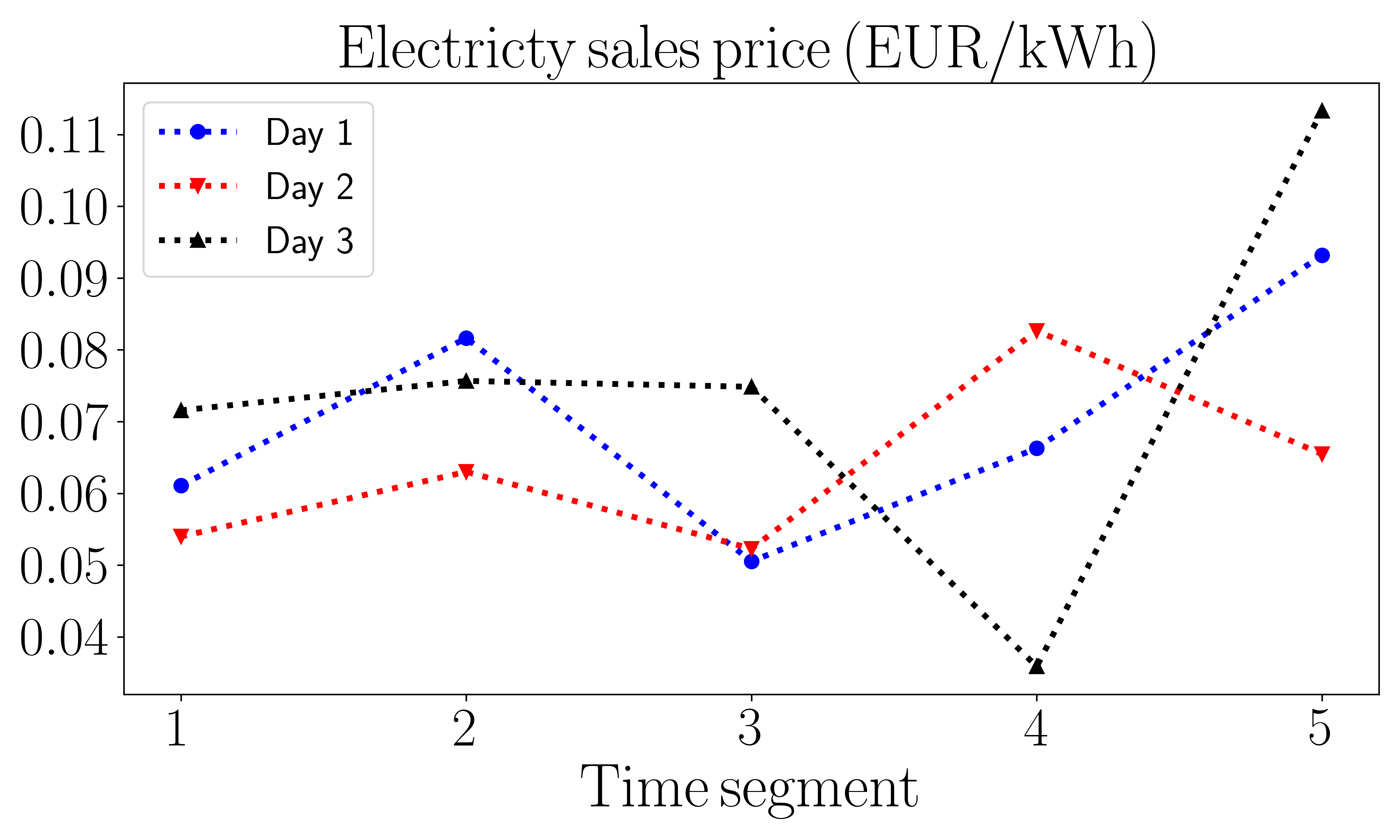}	
	\caption{The clustered electricity wholesale price with three representative days, each split into five time segments.}
	\label{fig:clustered_electricitySale}
\end{figure}

The time series were clustered into 3 representative days, while each day is split into 5 time segments and covers a $24\unit{\,\hour}$ period, as explained in section \ref{subsec:time-dependency and time series aggregation}. The choice regarding the number of days and time segments was made based on a design convergence study which can be found in \ref{app:temporalAnalysis}. The clustered time series of the wholesale electricity price is visualized in figure \ref{fig:clustered_electricitySale}, showing strong price variations within each day. For completeness, the clustered time series of the outside temperature and the total heat demand are visualized in figures \ref{fig:clustered_outsideTemperature} and \ref{fig:clustered_demand}, respectively.

\FloatBarrier

\subsection{Optimal storage sizing and operation}
\label{subsec:GlobalOptimum}  

The proposed optimization methodology is demonstrated in this section. The optimization decides jointly on the optimal storage size (its volume) together with the operation of the CHP unit (heat supply) and of the heat storage (charging and discharging) for each time segment. 

The results show that the optimal storage size is 1038 m$^3$. The storage integration reduces the levelized cost of heat (LCOH) by 16.5\% to $1.52\,\si{\sieuro cent / \kilo\watt\hour}$, down from $1.82\,\si{\sieuro cent / \kilo\watt\hour}$ for the reference case without storage, over $\npvNValue$ years, as shown in table \ref{tab:CostComparison}. The storage integration achieves a reduction in the heat production cost (CHP OPEX) of $1.07\,\si{\,\mega\sieuro}$, from $4.4\,\si{\,\mega\sieuro}$ down to $3.33\,\si{\,\mega\sieuro}$. The storage investment cost (Storage CAPEX) on the other hand amounts only to $323\,\si{\,\kilo\sieuro}$, resulting in a net financial benefit for the DHN operator of $730\,\si{\,\kilo\sieuro}$.

\begin{table}[t]
	\centering
	\caption{Economic comparison between the DHN with optimal storage integration and the reference operation without storage.}	
	\label{tab:CostComparison}
	\begin{tabularx}{\columnwidth}{ L{0.5} L{0.2} L{0.3} }
		& With \linebreak storage & Without storage\\
		\hline
		\hline		
		Pump OPEX $\pro$ $(\si{\,\mega\sieuro})$  & 0.0129 & 0.017 \\
		Pump OPEX $\stor$ $(\si{\,\mega\sieuro})$ & 0.0162 & - \\
		CHP OPEX $(\si{\,\mega\sieuro})$   & 3.33   & 4.4 \\
		Storage CAPEX $(\si{\,\mega\sieuro})$     & 0.323  & - \\
		Total cost $(\si{\,\mega\sieuro})$        & 3.68   & 4.41 \\
		\hline
		LCOH $(\si{\sieuro cent / \kilo\watt\hour})$                & 1.52   & 1.82\\
		\hline
		Storage size (m$^3$)                      & 1038   & - \\
	\end{tabularx}
\end{table} 

\begin{figure}[t]
	\centering
	%singlecolumn figure
	\includegraphics[width=\columnwidth]{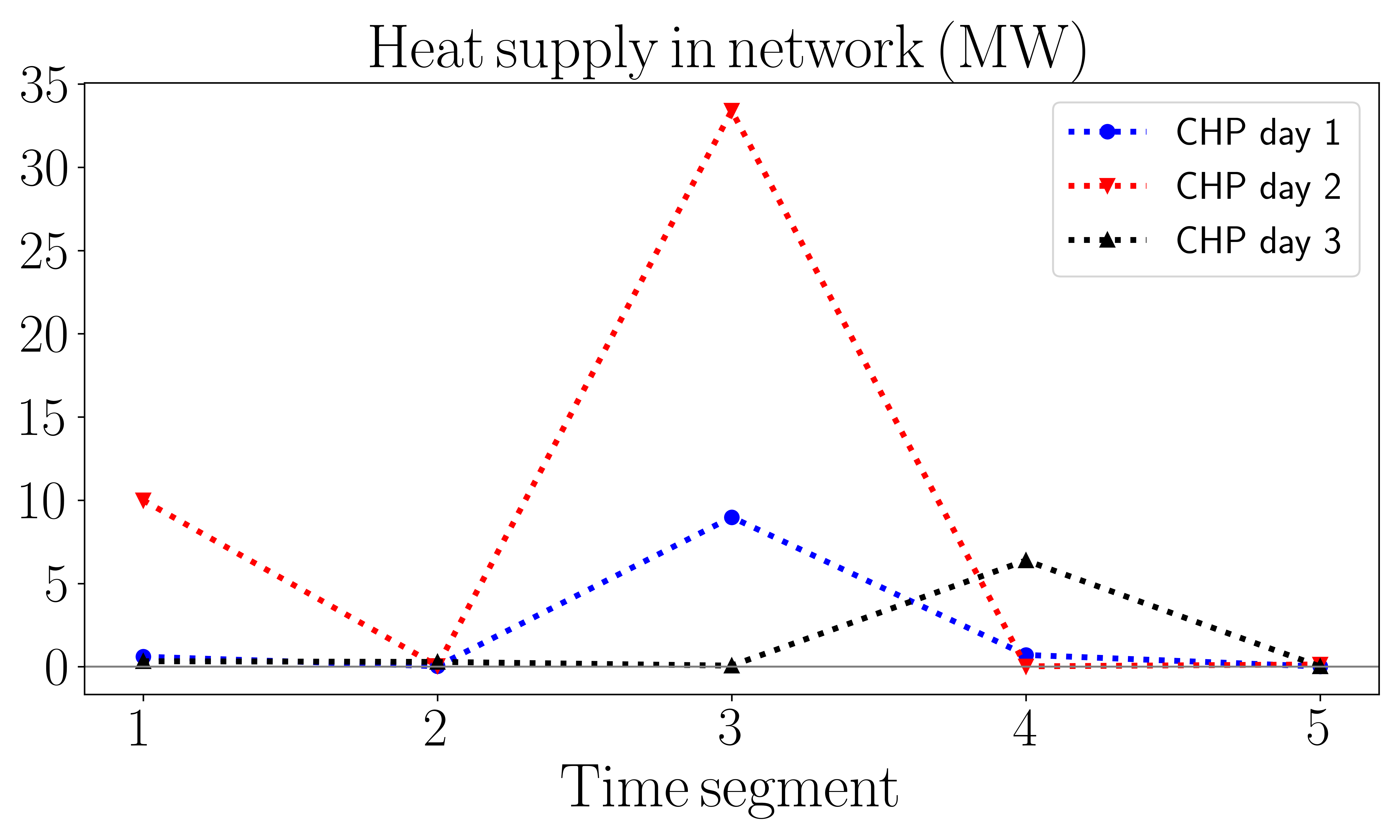}	
	\caption{Optimized heat supply to the DHN by the CHP unit for the case with storage.}
	\label{fig:Production_operation}
\end{figure}

\begin{figure}[t]
	\centering
	%singlecolumn figure
	\includegraphics[width=\columnwidth]{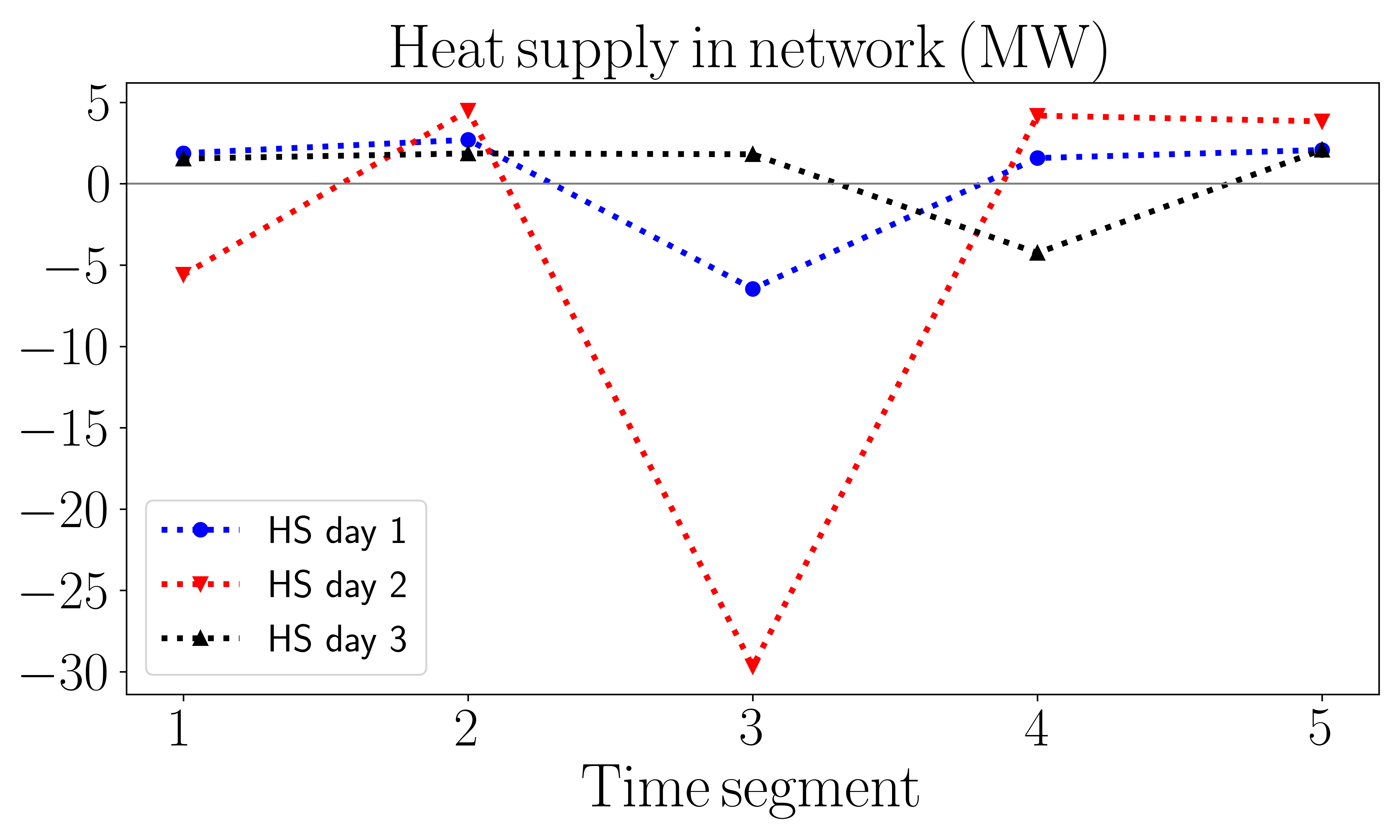}	
	\caption{Optimized heat supply to the DHN by the heat storage (HS). Note that negative values indicate a charging of the storage.}
	\label{fig:Storage_operation}
\end{figure}

The cost reduction in the CHP OPEX for the case with storage is achieved as the storage allows to shift the heat extraction from the CHP to moments of low electricity wholesale prices. This can be seen in figures \ref{fig:Production_operation} and \ref{fig:Storage_operation} where the operation of the CHP unit and the storage are shown, respectively. Looking for example at the representative day 3, it can be seen that the CHP unit supplies heat only during the time segment 4, where the electricity price is lowest (see figure \ref{fig:clustered_electricitySale}), while the heat storage provides the required heat to the DHN during the other time segments of that day and is charged only during time segment 4.

To further analyze the optimal storage sizing, the relation between the storage size and the costs are studied in detail through a parameter scan. Here, the storage size is fixed at values between 0 m$^3$ (no storage) up to 2500 m$^3$. For each size, the operation of the CHP and storage are optimized.

\begin{figure}[t]
	\centering
	%singlecolumn figure
	\includegraphics[width=\columnwidth]{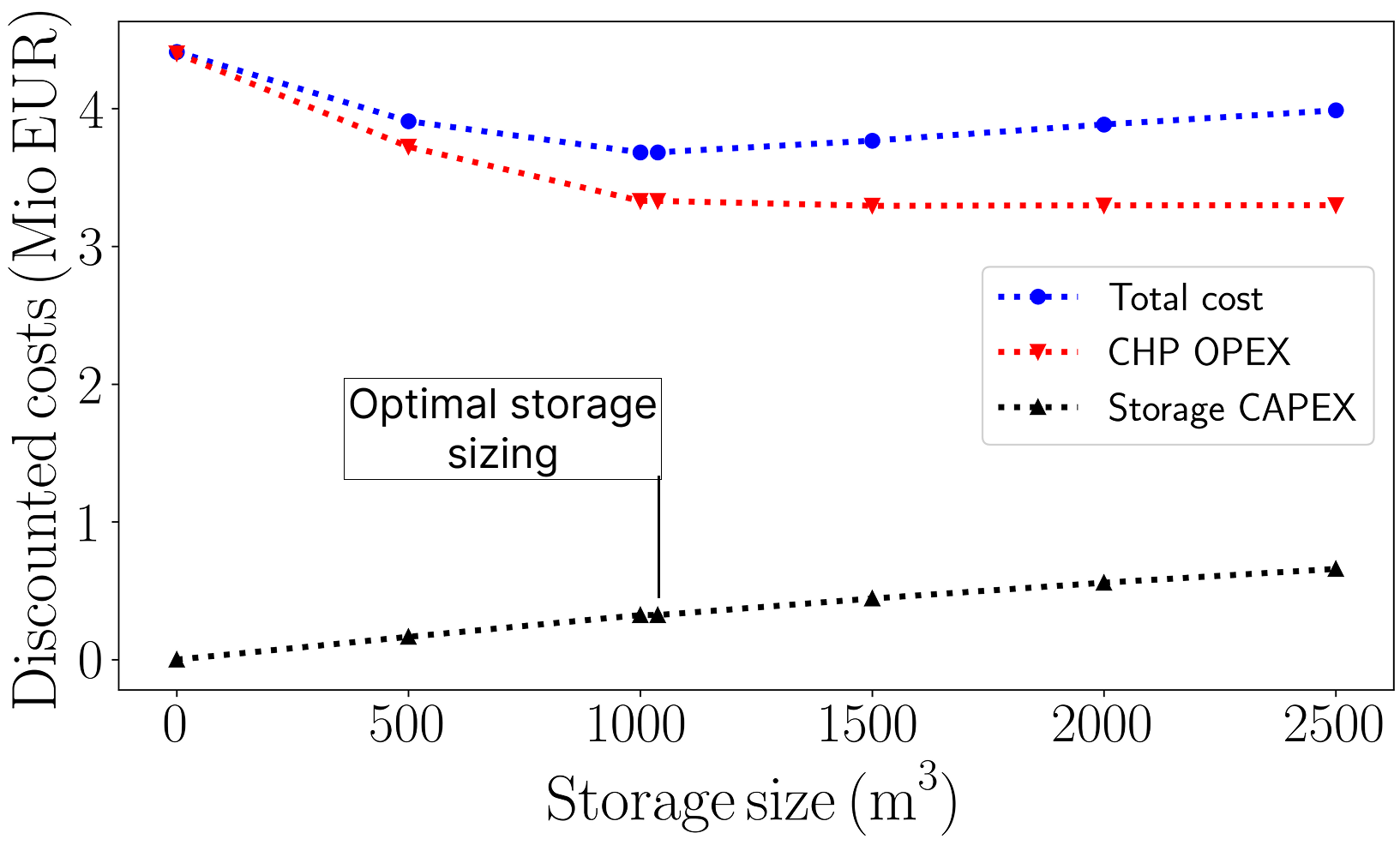}	
	\caption{Storage volume parameter scan for 3 representative days and 5 time segments, showing the impact of the storage size on the discounted cost.}
	\label{fig:volume_study}
\end{figure}

The results of this parameter scan, shown in figure \ref{fig:volume_study}, confirm that the methodology identifies the optimal storage size. They also reveal that every studied storage size from 500 m$^3$ up to 2500 m$^3$ achieves lower total cost than the reference case without storage, highlighting the strong value of adding a storage to a CHP-DHN. Nevertheless, the parameter scan also shows that an oversized storage of 2500 m$^3$ would cost the DHN operator additional $305\,\si{\,\kilo\sieuro}$ in comparison to the optimal storage size of 1038 m$^3$, underlining the value of optimal sizing.

Additionally, the profitability of the proposed storage investment was further assessed by reducing the variability\footnote{The difference between the price at $\dayVar,\timeSegVar$ and the average price of day $\dayVar$.} of the clustered electricity sales prices. The results show that even with a 50\% reduction in price variability, the investment remains profitable and reduces the LCOH by 4\% compared to not having a storage.

The solved optimization problem to obtain the optimal storage size contains 993 design variables and 979 technological constraints, and was solved on a laptop in less than 3 hours.

\subsection{Storage sizing with a simplified network model}
\label{subsec:simplfiedDHNmodel}

In the proposed physics-based optimization method for optimal heat storage sizing in DHNs, the entire heating network is modeled accurately and jointly with the heat storage and heat producer. As outlined in the literature review in section \ref{subsec:previous work on storage sizing}, prior storage sizing studies typically rely on simplified network representations, most commonly by aggregating the entire DHN into a single demand point. To assess the value of integrating a detailed DHN model into the storage optimization, the proposed approach is compared against this simplified aggregation strategy, which reduces the network to the configuration illustrated in figure \ref{fig:DHN_withStorage}, where the sum of the individual consumer heat demands is used for the demand of the aggregated consumer at each time segment.

The storage optimization based on the simplified approach yields a storage size of 771 m$^3$, which is 25.7\% smaller than the optimal storage size of 1038 m$^3$ obtained with the proposed method. This reduction stems from the fact that the pipe heat losses of the DHN are not captured by the simplified approach, leading to an underestimation of the required heat supply and consequently a smaller storage size. Network heat losses in DHNs can be significant, reaching up to 40\% of supplied heat in existing systems \cite{HeatLosses_iea2014}. For the DHN studied here, the annual network heat losses amount to 25.7\%\footnote{The physics-based DHN simulation, and thereby the calculation of the heat losses, was validated with detailed measurement data of flow rates and temperatures throughout the entire network and for a full-year operation, which was provided by the network operator.}, precisely the same percentage by which the simplified approach underestimates the storage size. In the simplified (aggregated) network model, pipe heat losses are also modeled, but they are practically zero due to the short distances between the producer, the storage, and the aggregated consumer. 

When the undersized storage of 771 m$^3$ is evaluated using the detailed DHN model while optimizing the producer and storage operation, the reduced storage capacity is found to limit the exploitation of varying electricity prices, increasing the total discounted cost by $132\,\si{\,\kilo\sieuro}$ or 3.6\% over $\npvNValue$ years compared to the optimal storage size of 1038 m$^3$.

While the simplified approach could potentially be improved for existing DHNs by incorporating measured (if available) or modeled heat losses into the aggregated demand, this approach remains limited in applicability. When either the storage location is subject to optimization, the DHN does not exist yet and needs to be designed as well, or when a storage is to be installed at the periphery of a DHN, detailed spatially-resolved knowledge of the network's physical state is essential to assess the interactions between the DHN and the storage. This information cannot be provided by simplified approaches. The proposed method, in contrast, retains full physical information of the network and, through the joint formulation of accurate DHN simulation and storage optimization, avoids time-consuming and potentially sub-optimal iterations between isolated simulation and optimization steps, making it a promising basis for more complex heat storage integration studies.

\FloatBarrier

\section{Discussion}
\label{sec:discussion} 
This paper presents an automated methodology based on mathematical optimization for cost-optimal heat-storage sizing in DHNs. The physics-based approach minimizes DHN operator costs by adding and optimally sizing a short-term storage unit, using a nonlinear flow and energy-transport model that captures network temperature and pressure drops and accounts for storage heat losses.

The methodology is applied to a 3rd-generation CHP-DHN in Oostende, Belgium, comprising 30 consumers. This exceeds the scale of existing physics-based storage sizing studies, which consider no more than 10 consumers and rely on manual choices or heuristic optimization rather than scalable mathematical optimization. Because of the underlying scalable optimization algorithms, the transformation of an intrinsically discrete problem into a continuous one, and the quasi-steady-state network model, this method is expected to scale well to larger DHNs with multiple potential storage sites.

For the considered case, the optimal storage integration reduces total costs (CAPEX and OPEX) by $730\si{\,\kilo\sieuro}$, or 16.5\%, to $3.68\si{\,\mega\sieuro}$, compared to $4.41\si{\,\mega\sieuro}$ without storage over $\npvNValue$ years. Savings arise from a reduction in CHP OPEX of $1.07\si{\,\mega\sieuro}$, from $4.4\si{\,\mega\sieuro}$ to $3.33\si{\,\mega\sieuro}$, as the storage shifts heat extraction to periods of low electricity wholesale prices. The storage CAPEX amounts to $323\si{\,\kilo\sieuro}$.

Additionally, the proposed methodology is compared to a simplified approach where the entire DHN is aggregated into a single demand point, as commonly used in previous storage sizing studies. The comparison reveals that the simplified approach underestimates the optimal storage size by 25.7\%, thereby limiting the exploitation of varying electricity prices, as pipe heat losses are not represented. Beyond pipe heat losses, the optimal storage size depends on the impact of CHP heat extraction on electricity production, temporal variation in electricity wholesale prices and heat demand, and storage investment cost. This highlights the advantage of an automated, physics-based optimization over simplified methods: it determines the cost-optimal storage size while explicitly accounting for heat losses in the DHN and storage, ensuring feasibility, accurate cost assessment, and optimal storage sizing. The optimization is computationally efficient, completing in under 3 hours on a standard laptop.

Future work should investigate additional applications of short-term storage in DHNs. These applications could include peak shaving of morning peaks after night setbacks and integrating renewable heat sources, such as solar thermal. In these cases, both storage size and location matter, as it affects costs, heat losses, and pipe sizing in new networks. The proposed optimization framework is well suited for such analyses because it models accurately DHN physics and storage–producer interactions in a nonlinear formulation. 

The 24-hour representative days used in this work limit the analysis of longer charging cycles. Considering multi-day or multi-week operation may change optimal storage sizes, as shown by \citet{SIFNAIOS2025}.

More detailed 1D storage models, such as multi-node or adaptive-grid formulations \cite{JODEIRI2024}, could better capture temperature profiles and improve assessments of charging strategies and storage sizing. The current 2-zone model may overestimate upper-layer temperatures because it neglects heat exchange between the hot and cold zones.

Operational limits on heat extraction should also be considered, since rapid changes of more than 30 MW within hours may not always be feasible and could affect storage economics.

Finally, uncertainties such as electricity price variability and DHN consumer expansion can influence optimal storage sizing. The sensitivity results in this work indicate that even a substantially oversized storage (e.g., 2500 m$^3$ vs. 1038 m$^3$) remains profitable. As larger storage volumes improve DHN flexibility \cite{GUELPA2019}, they may offer cost-effective robustness under uncertain future conditions, though further research is needed to confirm this.

\section{Conclusions and outlook}
\label{sec:conclusions}

This work presents a scalable, automated optimization methodology for cost-optimal integration of short-term heat storage in DHNs. By combining a nonlinear, physics-based thermal-hydraulic model with mathematical optimization, the method determines economically optimal storage sizes while ensuring feasible DHN and storage operation, as well as accurate cost assessment.

Applied to a real 3rd-generation CHP-DHN with 30 consumers, a scale that exceeds existing physics-based storage-sizing studies, the method reduces total 20-year costs by $730\si{\,\kilo\sieuro}$ (16.5\%), from $4.41\si{\,\mega\sieuro}$ to $3.68\si{\,\mega\sieuro}$. Savings result from a $1.07\si{\,\mega\sieuro}$ decrease in CHP operating cost by shifting heat extraction to low electricity price periods, while storage investment amounts to $323\si{\,\kilo\sieuro}$. 

A comparison to a commonly used simplified storage sizing method, which fails to identify the optimal storage size, demonstrates the advantage of the proposed holistic, physics-based optimization approach in ensuring feasible designs, accurate cost assessments, and optimal storage sizing. The results show that optimal storage sizing depends on the interaction between CHP heat extraction and electricity production, temporal price and demand variability, network heat losses, and storage investment cost. 

The developed framework lays the basis for future methodological extensions to applications such as peak shaving and renewable-heat integration, where both storage size and location matter. Further research should also investigate longer representative periods, more detailed storage models, operational limits on producer heat extraction, and uncertainties such as electricity price variability and network expansion.

\section*{CRediT authorship contribution statement}
\textbf{Martin Sollich}: Conceptualization, Data curation, Formal analysis, Funding acquisition, Investigation, Methodology, Software, Visualization, Writing – original draft.
\textbf{Maarten Blommaert}: Conceptualization,  Funding acquisition, Methodology, Supervision, Writing – review \& editing.

\section*{Declaration of competing interest}
The authors declare that they have no known competing financial interests or personal relationships that could have appeared to influence the reported work.

\section*{Data availability}
A data set containing the DHN structure and the time series used in the case studies of this paper is available at the following link: \url{https://doi.org/10.48804/OANORS}. The optimization results can be replicated using the methodology described in this paper.

\section*{Acknowledgments}
Martin Sollich is funded by the Research Foundation – Flanders (FWO) through the PhD fellowship strategic basic research with the file number 1SH7624N.

Maarten Blommaert and Martin Sollich have received funding from the KU Leuven, Belgium with the reference STG/21/016.

The authors would like to extend their gratitude to the operator of the studied DHN in Oostende, Beauvent cv, for their input to this work.

\appendix

\section{Key model and optimization parameters}\label{app:parameters}
\setcounter{figure}{0}
\setcounter{table}{0}
\setcounter{equation}{0}
\renewcommand{\thefigure}{\thesection\arabic{figure}}
\renewcommand{\thetable}{\thesection\arabic{table}}
\renewcommand{\theequation}{\thesection\arabic{equation}}

The key parameters used for the network model (simulation) and the storage optimization to obtain the results presented in section \ref{sec:case study} are given in table \ref{tab:properties}.

\newcommand{\CHPturbineEffValue}{0.2}
\newcommand{\pumpEffValue}{0.81}

\begin{table}[h]
	\centering
	\caption{Key parameters used for the network model (simulation) and the storage optimization.}	
	\label{tab:properties}
	\begin{tabularx}{\columnwidth}{llll}		
		Property & Value & Unit & Reference\\
		\hline

		$\CHPturbineEff$ & $\CHPturbineEffValue$ & $-$&DHN operator\\
		$\pumpEff$ & $\pumpEffValue$ & $-$&\\
		$\density$ & $\densityValue$ & $\unit{\kilogram\per\meter^3}$&\\
		
		$\spHeatCap$ & $\spHeatCapValue$ & $\unit{\joule\per(\kilogram\kelvin)}$&\\
		$\cPumpOPEXi$ & 0.1 & $\si{\sieuro / \kilo\watt\hour}$&DHN operator\\
		$\KOPEX$ & 8760 & $\unit{\hour\per\year}$ & \\
		$\maxPressure$ & 15 & $\unit{\bar}$ & \\
		$\Ustor$ & 0.12 & $\unit{\watt\per(\meter^2\kelvin)}$ & \cite{ROMANCHENKO2018}\\
		$\storMaxVolume$ & 6000 & $\unit{\meter^3}$&\\

	\end{tabularx}
\end{table}
\FloatBarrier

\section{Storage CAPEX fit}
\label{app:fitParameters}
\setcounter{figure}{0}
\setcounter{table}{0}
\setcounter{equation}{0}
\renewcommand{\thefigure}{\thesection\arabic{figure}}
\renewcommand{\thetable}{\thesection\arabic{table}}
\renewcommand{\theequation}{\thesection\arabic{equation}}

As a supplement to section \ref{sec:method} the details of the storage investment cost fit are provided in this section.

%% storage capex fit and parameters
\newcommand{\aStorCAPEX}{a_{C,\stor}}
\newcommand{\bStorCAPEX}{b_{C,\stor}}
\newcommand{\cStorCAPEX}{c_{C,\stor}}
\newcommand{\aStorCAPEXValue}{1.557e-06}
\newcommand{\bStorCAPEXValue}{-0.03848}
\newcommand{\cStorCAPEXValue}{349.9}

\begin{table}[h]
	\centering
	\caption{Parameters of the storage investment cost fit ($\storCostInv$) as shown in figure \ref{fig:storage_CAPEX_plot}. The polynomial fit reads $\storCostInv = \aStorCAPEX x^3 + \bStorCAPEX x^2 + \cStorCAPEX x$ with $x$ being the storage volume in $m^3$. The fit is based on data from academic literature and storage manufacturers \cite{LorenzBehaelterbau2014, Espagnet2016MasterThesis, HuchBehaelterbau2021}.}	
	\label{tab:storageCAPEXfit}
	\begin{tabularx}{\columnwidth}{XX}
		Coefficient & Value \\
		\hline
		$\aStorCAPEX$    & $\aStorCAPEXValue$ \\
		$\bStorCAPEX$    & $\bStorCAPEXValue$ \\
		$\cStorCAPEX$    & $\cStorCAPEXValue$ \\
	\end{tabularx}
\end{table}
\FloatBarrier

\section{Clustered parameters}
\label{app:timeSeries}
\setcounter{figure}{0}
\setcounter{table}{0}
\setcounter{equation}{0}
\renewcommand{\thefigure}{\thesection\arabic{figure}}
\renewcommand{\thetable}{\thesection\arabic{table}}
\renewcommand{\theequation}{\thesection\arabic{equation}}

\begin{figure}[h]
	\centering
	%singlecolumn figure
	\includegraphics[width=\columnwidth]{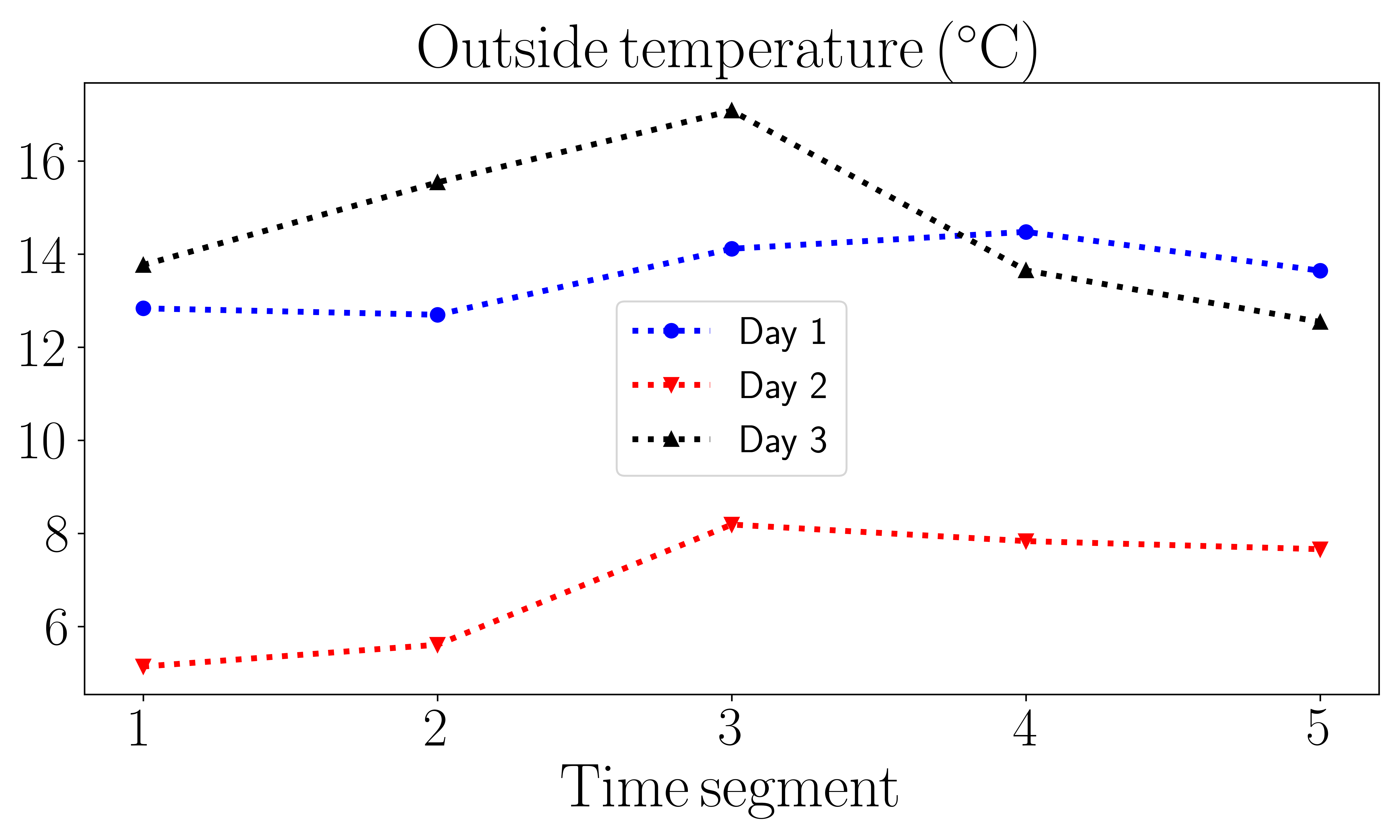}	
	\caption{The clustered outside temperature with three representative days, each with five time segments.}
	\label{fig:clustered_outsideTemperature}
\end{figure}

\begin{figure}[h]
	\centering
	%singlecolumn figure
	\includegraphics[width=\columnwidth]{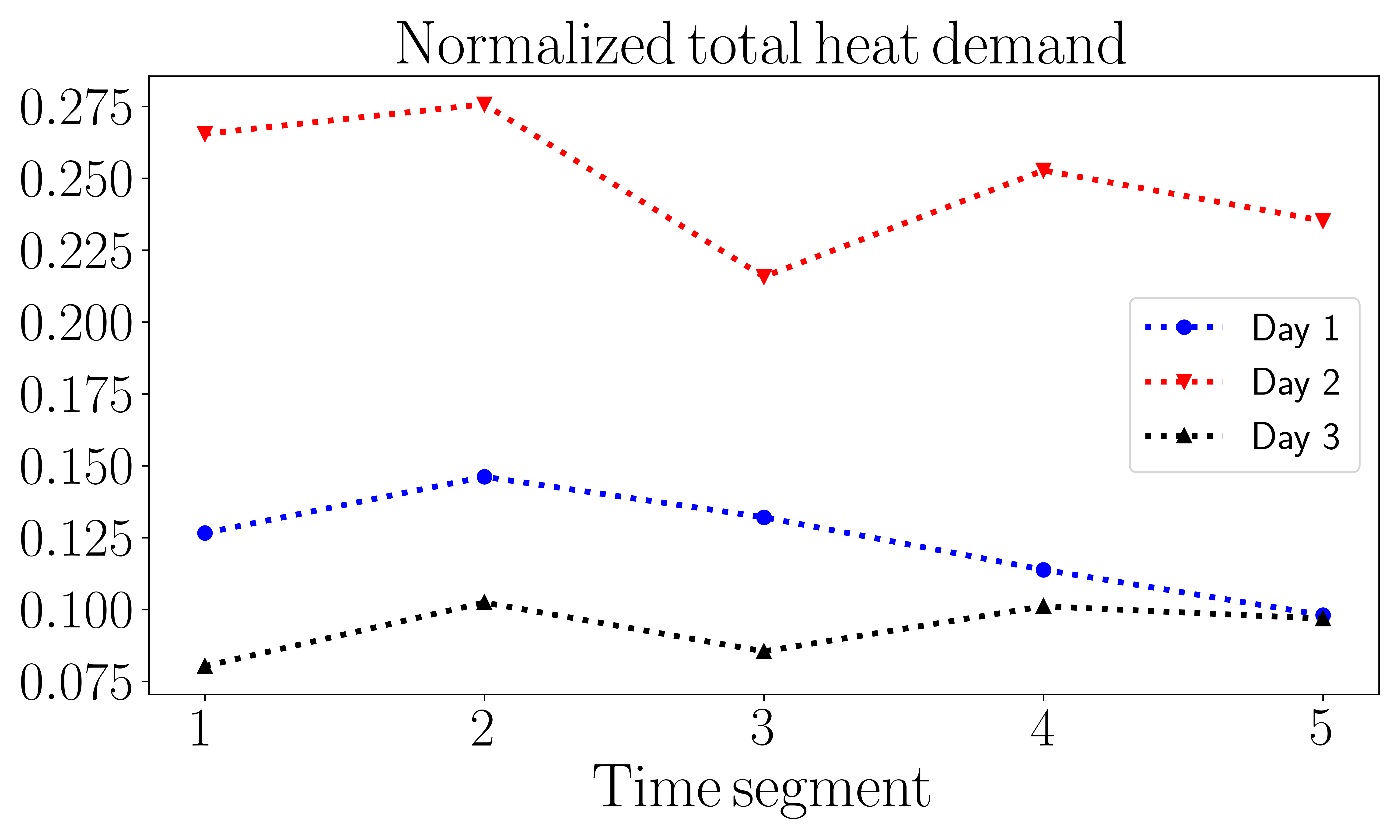}	
	\caption{The clustered total heat demand of the DHN with three representative days, each with five time segments.}
	\label{fig:clustered_demand}
\end{figure}

As a supplement to section \ref{sec:case study}, the clustered time series of the outside temperature and of the total normalized DHN heat demand are provided in figures \ref{fig:clustered_outsideTemperature} and \ref{fig:clustered_demand}, respectively.

\FloatBarrier

\section{Temporal clustering analysis}
\label{app:temporalAnalysis}  
\setcounter{figure}{0}
\setcounter{table}{0}
\setcounter{equation}{0}
\renewcommand{\thefigure}{\thesection\arabic{figure}}
\renewcommand{\thetable}{\thesection\arabic{table}}
\renewcommand{\theequation}{\thesection\arabic{equation}}

\begin{figure}[h]
	\centering
	%singlecolumn figure
	\includegraphics[width=\columnwidth]{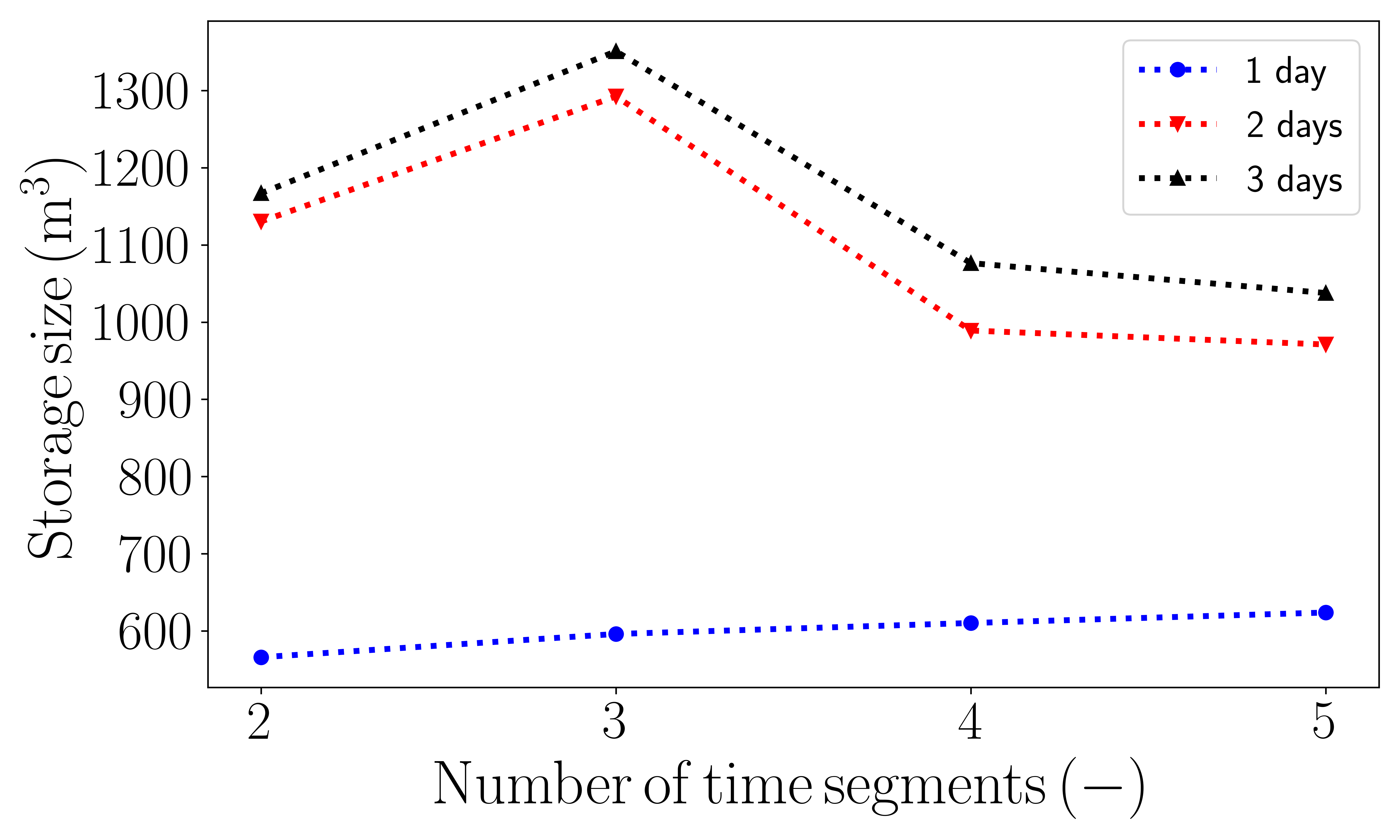}	
	\caption{Impact of the temporal resolution of the clustering, i.e., number of representative days and time segments, on the storage size.}
	\label{fig:timeStep_study_storageSize}
\end{figure}

To select the number of representative days and time segments per day, a sensitivity study was performed to assess their impact on the optimal storage size. The temporal resolution was gradually increased until design convergence was reached, i.e. the storage size changed only marginally, similar to \cite{CERUTI2025}. As shown in figure \ref{fig:timeStep_study_storageSize}, increasing from 2 to 3 representative days and from 4 to 5 time segments per day alters the optimal storage size by less than 100 m$^3$, respectively. 

\begin{figure}[h]
	\centering
	%singlecolumn figure
	\includegraphics[width=\columnwidth]{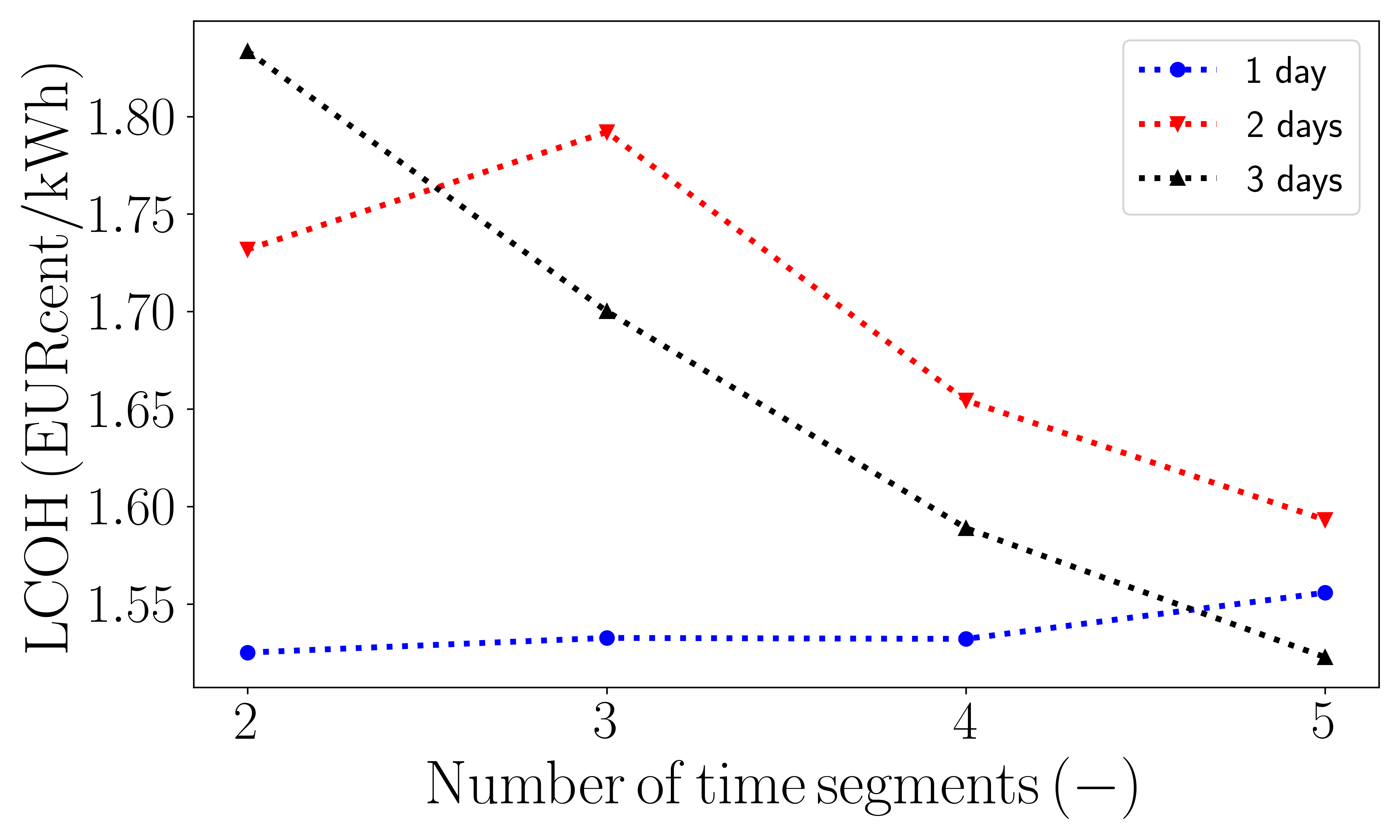}	
	\caption{Impact of the temporal resolution of the clustering, i.e., number of representative days and time segments, on the LCOH.}
	\label{fig:timeStep_study_LCOH}
\end{figure}

The effect on the LCOH, shown in \ref{fig:timeStep_study_LCOH}, is slightly more pronounced. Finer temporal resolution yields more extreme electricity-price values, enabling the optimization to charge the storage during very low-price segments and achieve additional savings. However, daily heat demand quantities change little, so the storage volume remains largely unaffected. The LCOH differences between 2 and 3 representative days (at 5 segments) or between 4 and 5 segments (at 3 days) remain below $0.07\,\si{\sieuro cent / \kilo\watt\hour}$, i.e. below 5\% of the LCOH for 3 days and 5 segments.

Since storage sizing is the focus of this work, 3 representative days and 5 time segments were selected. Moreover, \citet{CERUTI2025} found that most accuracy gains, while maintaining low computational cost, are achieved when considering at least 2 representative days and 5 or more segments, consistent with the findings of this study.

\FloatBarrier

\bibliographystyle{elsarticle-num-names}
\bibliography{library}% common bib file

\end{document}